\documentclass[11pt]{amsart}
\usepackage{amssymb}
\usepackage{amsmath}
\usepackage{fancyhdr}

\title[Hamilton decompositions of regular expanders]{Hamilton decompositions of regular expanders: applications
}

\date{\today}
\author{Daniela K\"uhn and Deryk Osthus}
\thanks{The authors were supported by the EPSRC, grant no.~EP/J008087/1. D. K\"uhn was also supported by the ERC, grant no.~258345.} 

\newtheorem{firstthm}{Proposition}[section]
\newtheorem{theorem}[firstthm]{Theorem}
\newtheorem{prop}[firstthm]{Proposition}
\newtheorem{lemma}[firstthm]{Lemma}
\newtheorem{cor}[firstthm]{Corollary}

\newtheorem{conj}[firstthm]{Conjecture}

\def\noproof{{\unskip\nobreak\hfill\penalty50\hskip2em\hbox{}\nobreak\hfill%
       $\square$\parfillskip=0pt\finalhyphendemerits=0\par}\goodbreak}
\def\endproof{\noproof\bigskip}
\newdimen\margin   % needed for macros \textdisplay & \ltextdisplay
\def\textno#1&#2\par{%
   \margin=\hsize
   \advance\margin by -4\parindent
          \setbox1=\hbox{\sl#1}%
   \ifdim\wd1 < \margin
      $$\box1\eqno#2$$%
   \else
      \bigbreak
      \hbox to \hsize{\indent$\vcenter{\advance\hsize by -3\parindent
      \it\noindent#1}\hfil#2$}%
      \bigbreak
   \fi}
\def\proof{\removelastskip\penalty55\medskip\noindent{\bf Proof. }}

\begin{document}

\def\COMMENT#1{}
\def\TASK#1{}

\def\eps{{\varepsilon}}
\newcommand{\ex}{\mathbb{E}}
\newcommand{\pr}{\mathbb{P}}
\newcommand{\cB}{\mathcal{B}}
\newcommand{\cS}{\mathcal{S}}
\newcommand{\cF}{\mathcal{F}}
\newcommand{\cC}{\mathcal{C}}
\newcommand{\cP}{\mathcal{P}}
\newcommand{\cQ}{\mathcal{Q}}
\newcommand{\cR}{\mathcal{R}}
\newcommand{\cK}{\mathcal{K}}
\newcommand{\cD}{\mathcal{D}}
\newcommand{\cI}{\mathcal{I}}
\newcommand{\cV}{\mathcal{V}}
\newcommand{\1}{{\bf 1}_{n\not\equiv \delta}}
\newcommand{\eul}{{\rm e}}

\begin{abstract}  \noindent
In a recent paper, we showed that every sufficiently large regular digraph $G$ on $n$ vertices whose degree is linear in $n$ and which is
a robust outexpander has a decomposition into edge-disjoint Hamilton cycles. 
The main consequence of this theorem is that every regular tournament on $n$ vertices can be decomposed
into $(n-1)/2$ edge-disjoint Hamilton cycles, whenever $n$ is sufficiently large. This verified a conjecture of Kelly from 1968.
In this paper, we derive a number of further consequences of our result on robust outexpanders, the main ones are the following:
\begin{itemize}
\item[(i)] an undirected analogue of our result on robust outexpanders;
\item[(ii)] best possible bounds on the size of an optimal packing of edge-disjoint Hamilton cycles in a graph of minimum degree $\delta$
for a large range of values for $\delta$.
\item[(iii)] a similar result for digraphs of given minimum semidegree;
\item[(iv)] an approximate version of a conjecture of Nash-Williams on Hamilton decompositions of dense regular graphs;
%\item[(v)] the observation that dense quasi-random graphs are robust outexpanders;
\item[(v)] a verification of the `very dense' case of a conjecture of Frieze and Krivelevich on packing edge-disjoint Hamilton cycles in random graphs;
\item[(vi)]  a proof of a conjecture of Erd\H{o}s on the size of an optimal packing of edge-disjoint Hamilton cycles in a random tournament.
\end{itemize}
\end{abstract}
\maketitle
\begin{center}
{\large \emph{accepted for publication in J. Combinatorial Theory B}}
\end{center}
\section{Introduction}\label{intro}

\subsection{Robust expanders}
A graph or digraph $G$ has a Hamilton decomposition if it contains a set of edge-disjoint Hamilton cycles which together
cover all the edges of $G$. Despite the fact that the study of Hamilton decompositions has a long history,
general results on Hamilton decompositions have been rare so far.
The first result in this direction is Walecki's construction of a Hamilton decomposition for the complete graph on an odd number of vertices.
Tillson~\cite{till} proved an analogue of this for complete digraphs.
In 1968, Kelly (see e.g.~\cite{bondy,KOsurvey,moon}) conjectured an analogue of this for tournaments, namely that every regular tournament has a Hamilton decomposition.
Note that in a digraph we allow up to two edges between any pair of vertices -- at most one in each direction. So a digraph might contain cycles of length two, whereas in an oriented graph, 
we only allow at most one edge between any pair of vertices. A tournament $T$ is an orientation of a complete (undirected) graph.
$T$ is regular if the outdegree of every vertex equals its indegree.

Very recently, we were able to prove Kelly's conjecture for all large tournaments~\cite{monster} 
(that paper also includes a more thorough discussion of partial results towards this conjecture).
Our proof led us to a far more general decomposition theorem, which involves the concept of robust outexpanders.
This concept was introduced by K\"uhn, Osthus and Treglown~\cite{KOTchvatal}, and was already used implicitly in~\cite{keevashko,kelly}.
Roughly speaking, a digraph is a robust outexpander if for every set $S$ which is not too small and not too large, 
its `robust' outneighbourhood is at least a little larger than $S$.
In~\cite{monster} we proved that every dense regular robust outexpander has a Hamilton decomposition.
Since every regular tournament is easily seen to be a robust outexpander, our result on Kelly's conjecture is a very special case of this.

More precisely, let $0<\nu\le  \tau<1$. Given any digraph~$G$ on $n$ vertices and
$S\subseteq V(G)$, the \emph{$\nu$-robust outneighbourhood~$RN^+_{\nu,G}(S)$ of~$S$}
is the set of all those vertices~$x$ of~$G$ which have at least $\nu n$ inneighbours
in~$S$. $G$ is called a \emph{robust $(\nu,\tau)$-outexpander}
if 
$$
|RN^+_{\nu,G}(S)|\ge |S|+\nu n \ \mbox{ for all } \ 
S\subseteq V(G) \ \mbox{ with } \ \tau n\le |S|\le  (1-\tau)n.
$$
We can now state the main result of~\cite{monster}, which guarantees a Hamilton decomposition in every regular robust outexpander of linear degree.
\begin{theorem} \label{decomp}
For every $ \alpha >0$ there exists $\tau>0$ such that for every $\nu >0$ there exists $n_0=n_0 (\alpha,\nu,\tau)$ 
for which the following holds. Suppose that
\begin{itemize}
\item[{\rm (i)}] $G$ is an $r$-regular digraph on $n \ge n_0$ vertices, where $r\ge \alpha n$;
\item[{\rm (ii)}] $G$ is a robust $(\nu,\tau)$-outexpander.
\end{itemize}
Then $G$ has a Hamilton decomposition.
Moreover, this decomposition can be found in time polynomial in $n$.
\end{theorem}
The proof of Theorem~\ref{decomp} introduces a new method for decomposing a graph into (Hamilton) cycles, which we believe will 
have further applications. 
As a tool, it uses (amongst others) a result of Osthus and Staden~\cite{OS} which states that any $G$ as in Theorem~\ref{decomp} has an approximate Hamilton decomposition
-- i.e.~a set of edge-disjoint Hamilton cycles covering almost all edges of $G$
(\cite{OS} generalizes a corresponding result in~\cite{KOTkelly} on dense regular oriented graphs).
Given such an approximate decomposition of $G$, it is obviously not always possible to 
extend this into a Hamilton decomposition. 
The key idea in the proof of Theorem~\ref{decomp} is to overcome this obstacle by  finding a sparse `robustly decomposable' spanning subdigraph $H^*$ of $G$, 
i.e.~$H^*$ has a Hamilton decomposition even if we add the edges of any very sparse spanning regular subdigraph $H$ of $G\setminus E(H^*)$ to $H^*$.
Now one can apply the result of~\cite{OS} to obtain an approximate decomposition of  $G\setminus E(H^*)$, which leaves a very sparse uncovered leftover $H$.
The choice of $H^*$ implies that  $H^* \cup H$ has a Hamilton decomposition. Altogether, this gives the required Hamilton decomposition of $G$.

Some applications of Theorem~\ref{decomp} were already derived and described in~\cite{monster}
(namely Theorems~\ref{regdigraph} and~\ref{orientcor} below). As also discussed in~\cite{monster}, Theorem~\ref{regdigraph} in turn has a 
surprising application to finding tours for the Asymmetric Travelling Salesman Problem with large `domination ratio', which (together with a result 
by Gutin and Yeo~\cite{GutinYeo}) solved a problem 
raised e.g.~by  Glover and Punnen~\cite{GP} as well as Alon, Gutin and Krivelevich~\cite{AGK}. 
Since robust expansion is a property shared by several important classes of graphs, it turns out that Theorem~\ref{decomp} has numerous further applications, which 
we discuss and derive in this paper.

For applications to undirected graphs, it is obviously helpful to have an undirected version of Theorem~\ref{decomp}.
For this, we introduce the analogue of robust outexpanders for undirected graphs. Let $0<\nu\le  \tau<1$. Given any
(undirected) graph~$G$ on $n$ vertices and
$S\subseteq V(G)$, the \emph{$\nu$-robust neighbourhood~$RN_{\nu,G}(S)$ of~$S$}
is the set of all those vertices~$x$ of~$G$ which have at least $\nu n$ neighbours
in~$S$. $G$ is called a \emph{robust $(\nu,\tau)$-expander}
if 
$$
|RN_{\nu,G}(S)|\ge |S|+\nu n  \ \mbox{ for all } \
S\subseteq V(G) \ \mbox { with } \ \tau n\le |S|\le  (1-\tau)n.
$$
With this notion, we can derive the following analogue for (undirected) robust expanders from Theorem~\ref{decomp}.
\begin{theorem}\label{undir_decomp}
For every $ \alpha >0$ there exists $\tau>0$ such that for every $\nu>0$ there exists $n_0=n_0 (\alpha,\nu,\tau)$ 
for which the following holds. 
Suppose that
\begin{itemize}
\item[{\rm (i)}] $G$ is an $r$-regular graph on $n \ge n_0$ vertices, where $r\ge \alpha n$ is even;
\item[{\rm (ii)}] $G$ is a robust $(\nu,\tau)$-expander.
\end{itemize}
Then $G$ has a Hamilton decomposition. Moreover, this decomposition can be found in time polynomial in $n$.
\end{theorem}
One can easily derive a version of the above result for odd values of~$r$:
suppose that $G$ satisfies all conditions of Theorem~\ref{undir_decomp} except that $r$ is odd.
Then we can obtain a decomposition of~$G$ into one perfect matching and a set of Hamilton cycles.
Indeed, the conditions ensure that $G$ contains a perfect matching~$M$ (e.g.~via Tutte's theorem).
Removing~$M$ from $G$ leaves a graph $G-M$ which satisfies the conditions of Theorem~\ref{undir_decomp}
(with slightly worse parameters), so $G-M$ has a Hamilton decomposition.

\subsection{Edge-disjoint Hamilton cycles in dense  graphs}
Nash-Williams~\cite{NW71} conjectured that every sufficiently dense even-regular graph has a Hamilton decomposition.
More precisely, he conjectured that such a decomposition exists in every $2d$-regular graph with at most $4d+1$ vertices.
This is sometimes referred to as the `Hamilton decomposition conjecture'.
The following result implies an approximate version of this conjecture.

\begin{theorem}\label{reggraph}
For every $\eps>0$ there exists $n_0$ such that every $r$-regular graph $G$ on $n\ge n_0$ vertices, where $r\ge (1/2+\eps)n$
is even, has a Hamilton decomposition.
\end{theorem} 
Theorem~\ref{reggraph} is an (almost) immediate consequence of Theorem~\ref{undir_decomp}.
After the current paper was completed, the conjecture of Nash-Williams was verified for all large~$n$ in a sequence of papers
by Csaba, K\"uhn, Lo, Osthus and Treglown~\cite{paper2,paper3,paper1,paper4}
(these papers also contain a proof of the related `$1$-factorization conjecture').
The proof is far more difficult than that of Theorem~\ref{reggraph}, but once again Theorem~\ref{undir_decomp} plays a crucial role.
Previous partial results related to this problem
were proved e.g.~Jackson~\cite{jackson}, Perkovic and Reed~\cite{PRedges}, as well as Christofides, K\"uhn and Osthus~\cite{CKO}.

The paper~\cite{CKO} also  asymptotically determined the number of edge-disjoint Hamilton cycles which one can guarantee in
a graph with given minimum degree~$\delta$. (The case $\delta=n/2$ of this question was originally raised by 
Nash-Williams~\cite{Nash-Williams70,NW71}.)
More precisely, given a graph $G$, let ${\rm ham}(G)$ denote the maximum number of edge-disjoint Hamilton cycles in $G$
and let ${\rm reg}_{\rm even}(G)$ denote the largest degree of an even-regular spanning subgraph of $G$.
So ${\rm ham}(G)\le {\rm reg}_{\rm even}(G)/2$. Moreover, given $\delta,n\in\mathbb{N}$ with $\delta<n$, let ${\rm ham}(n,\delta)$
denote the minimum of ${\rm ham}(G)$ over all graphs on $n$ vertices with minimum degree~$\delta$.
Similarly, let ${\rm reg}_{\rm even}(n,\delta)$
denote the minimum of ${\rm reg}_{\rm even}(G)$ over all graphs on $n$ vertices with minimum degree~$\delta$.
Thus 
\begin{equation} \label{regham}
{\rm ham}(n,\delta)\le {\rm reg}_{\rm even}(n,\delta)/2
\end{equation} 
and
${\rm ham}(n,\delta)=0$ whenever $\delta<n/2$.
The following result was proved in~\cite{CKO}. It determines ${\rm reg}_{\rm even}(n,\delta)$ almost exactly and
gives a lower bound on ${\rm ham}(n,\delta)$ which approximately matches the trivial upper bound in~(\ref{regham}).
The lower bound in~(i) is Theorem~12(i) in~\cite{CKO} and
the upper bound in~(i) was proved in Section~2 of~\cite{CKO} (see also~\cite{Hartkefactors}).
\begin{theorem}\label{CKOthm}
Suppose that $n,\delta\in\mathbb{N}$ and $n/2< \delta<n$.
Let 
$$g(n,\delta):=\frac{\delta+\sqrt{n(2\delta-n)}}{2}.$$
Let $g_{\rm even}(n,\delta)$ be the largest even number which is at most $g(n,\delta)$ and let $g'_{\rm even}(n,\delta)$ be the largest even
number which is at most $g(n,\delta)+1$.
Then 
\begin{itemize}
\item[(i)]
%\begin{equation} \label{neu}
$g_{\rm even}(n,\delta)\le {\rm reg}_{\rm even}(n,\delta)\le g'_{\rm even}(n,\delta).$
%\end{equation}
\item[(ii)] For every $\eps>0$ there exists $n_0$ such that whenever $n\ge n_0$ and $(1/2+\eps)n< \delta<n$ then
$$ g_{\rm even}(n,\delta)-\eps n\le 2\cdot {\rm ham}(n,\delta)\le {\rm reg}_{\rm even}(n,\delta)\le g'_{\rm even}(n,\delta).$$
\end{itemize}
\end{theorem}
The error bound in (ii) was subsequently improved by Hartke and Seacrest~\cite{HartkeHCs}.
%Nash-Williams~\cite{Nash-Williams70} had originally conjectured that ${\rm ham}(n, \lceil n/2 \rceil)=\lfloor \delta/2 \rfloor$.
%However, Babai disproved this by observing that  ${\rm reg}_{\rm even}(n,\lceil n/2\rceil)$ is close to $n/4$.
Note that if $\delta$ is close to $n/2$, then $g(n,\delta)$ is close to $n/4$.
In particular, Theorem~\ref{CKOthm} implies an approximate solution to the problem raised by Nash-Williams.
The upper bounds in Theorem~\ref{CKOthm} 
(as well as those in Theorem~\ref{mindigraph} below) are based on a generalization of a construction of Babai.

Theorem~\ref{reggraph} together with Theorem~\ref{CKOthm}(i) imply the following precise version
of Theorem~\ref{CKOthm}(ii) for graphs of sufficiently large minimum degree.
Surprisingly, it turns out that the trivial bound~(\ref{regham}) holds with equality in this case.
Note that there are many pairs $n,\delta$ for which we have $g_{\rm even}(n,\delta) =g'_{\rm even}(n,\delta)$,
and thus Theorem~\ref{mingraph} even gives the  exact numerical value for ${\rm ham}(n,\delta)$ in these cases.
\begin{theorem}\label{mingraph}
For every $\eps>0$ there exists $n_0$ such that the following holds for all $n\ge n_0$. 
\begin{itemize}
\item[(i)] Every graph $G$ on $n$ vertices with $\delta(G) \ge (2-\sqrt{2}+\eps)n$ satisfies ${\rm ham}(G)= {\rm reg}_{\rm even}(G)/2$. 
\item[(ii)] If $(2-\sqrt{2}+\eps)n< \delta< n$, then
$$g_{\rm even}(n,\delta)\le 2\cdot {\rm ham}(n,\delta)= {\rm reg}_{\rm even}(n,\delta)\le g'_{\rm even}(n,\delta),$$
where $g_{\rm even}(n,\delta)$ and $g'_{\rm even}(n,\delta)$ are as defined in Theorem~\ref{CKOthm}.
\end{itemize}
\end{theorem}
Note that (i) immediately implies (ii).
In~\cite{paper2,paper3,KLOmindeg,paper1,paper4}, Theorem~\ref{mingraph}(ii) is extended to cover the entire range when $\delta \ge n/2$.
The proof relies on Theorem~\ref{mingraph} for the case when $\delta$ is much larger than $n/2$ and again relies on Theorem~\ref{undir_decomp}.
The paper by K\"uhn, Lapinskas and Osthus~\cite{KLOmindeg} covers the case when $G$ is far from extremal and
the sequence~\cite{paper2,paper3,paper1,paper4} mentioned above covers the `extremal cases'.
Together, these results imply an exact solution to the problem of Nash-Williams mentioned earlier.

A challenging open problem would be to extend the stronger assertion~(i) to the entire range $\delta(G) \ge n/2$.
An approximate result towards this was recently proved by Ferber, Krivelevich and Sudakov~\cite{FKS}.

\subsection{Edge-disjoint Hamilton cycles in dense digraphs}

We can use Theorem~\ref{decomp} to prove an analogue of Theorems~\ref{CKOthm} and~\ref{mingraph} for digraphs.
Before we can state it, we need to introduce the following notation.
Given a digraph $G$, we write $\delta^+(G)$ for its minimum outdegree and $\delta^-(G)$ for its minimum indegree.
The \emph{minimum semidegree} $\delta^0(G)$ of $G$ is the minimum of $\delta^+(G)$ and $\delta^-(G)$.

Similarly as for undirected graphs, given a digraph $G$, let ${\rm ham}(G)$ denote the maximum number of edge-disjoint Hamilton cycles in $G$
and let ${\rm reg}(G)$ denote the largest degree of a regular spanning subdigraph of $G$.
So ${\rm ham}(G)\le {\rm reg}(G)$. Moreover, given $\delta,n\in\mathbb{N}$ with $\delta<n$, let ${\rm ham}^{\rm dir}(n,\delta)$
denote the minimum of ${\rm ham}(G)$ over all digraphs on $n$ vertices with minimum semidegree $\delta$.
Similarly, let ${\rm reg}^{\rm dir}(n,\delta)$
denote the minimum of ${\rm reg}(G)$ over all digraphs on $n$ vertices with minimum semidegree $\delta$.
Thus ${\rm ham}^{\rm dir}(n,\delta)\le {\rm reg}^{\rm dir}(n,\delta)$ and ${\rm ham}^{\rm dir}(n,\delta)=0$ whenever
$\delta <n/2$.

\begin{theorem}\label{mindigraph}
Given
$\delta, n\in \mathbb{N}$ with $n/2\le \delta<n$, let
$$
f(n,\delta):=\left\lfloor\frac{\delta+\sqrt{n(2\delta-n)+\1}}{2}\right\rfloor, \ \ \ \text{where}  \ \ \
{\bf 1}_{n\not\equiv \delta}:=\begin{cases} 0 & \mbox{if } n\equiv \delta \mod 2\\
1 & \mbox{if } n\not \equiv \delta \mod 2.\end{cases}
$$
For every $\eps>0$ there exists $n_0$ such that the following statements hold for all $n\ge n_0$. 
\begin{itemize}
\item[{\rm (i)}] If $(1/2+\eps)n< \delta< n$ then
$$ f(n,\delta)-\eps n\le {\rm ham}^{\rm dir}(n,\delta)\le {\rm reg}^{\rm dir}(n,\delta)=f(n,\delta).$$
\item [{\rm (ii)}] Every digraph $G$ on $n$ vertices with $\delta^0(G) \ge (2-\sqrt{2}+\eps)n$ satisfies ${\rm ham}(G)= {\rm reg}(G)$. 
In particular, if $(2-\sqrt{2}+\eps)n< \delta<n$ then 
$${\rm ham}^{\rm dir}(n,\delta)={\rm reg}^{\rm dir}(n,\delta)=f(n,\delta).$$
\end{itemize}
\end{theorem}  
We conjecture that one can extend Theorem~\ref{mindigraph} to show that 
${\rm ham}^{\rm dir}(n,\delta)={\rm reg}^{\rm dir}(n,\delta)$ for any $\delta \ge n/2$.
This would completely determine the values of ${\rm ham}^{\rm dir}(n,\delta)$, as in Lemma~\ref{factor} we show that
${\rm reg}^{\rm dir}(n,\delta)=f(n,\delta)$ for any $\delta \ge n/2$.

We now state two consequences of Theorem~\ref{decomp} which we already observed in~\cite{monster} and which we will use in this paper.
The first is a version of Theorem~\ref{reggraph} for digraphs.
We will use it in our proof of Theorem~\ref{mindigraph}.
\begin{theorem}\label{regdigraph}
For every $\eps>0$ there exists $n_0$ such that every $r$-regular digraph $G$ on $n\ge n_0$ vertices, where $r\ge (1/2+\eps)n$,
has a Hamilton decomposition.
\end{theorem} 
Theorem~\ref{regdigraph} follows from Theorem~\ref{decomp} since
one can easily verify that regular digraphs as above are robust outexpanders (see Lemma~13.2~in~\cite{monster}).

The second result states that every sufficiently dense regular oriented graph has a Hamilton decomposition.
We will use it in our proof of Theorem~\ref{erdos} below.
\begin{theorem}\label{orientcor}
For every $\eps>0$ there exists $n_0$ such that every $r$-regular oriented graph $G$ on $n\ge n_0$
vertices with $r\ge 3n/8+\eps n$ has a Hamilton decomposition.
\end{theorem}
Note that this implies Kelly's conjecture for large regular tournaments.
As before, one can easily verify that such oriented graphs are robust outexpanders
(see Lemma~13.1~in~\cite{monster}).

\subsection{Random graphs and random tournaments}
Erd\H{o}s conjectured a probabilistic version of Kelly's conjecture (see~\cite{thomassen79}), namely that asymptotically almost surely a random tournament~$T$
should contain  $\delta^0(T)$ edge-disjoint Hamilton cycles. Here we say that a property holds asymptotically almost surely (a.a.s.) if
it holds with probability tending to~$1$ as $n$ tends to infinity.
Note that trivially, any digraph $G$ has at most $\delta^0(G)$ edge-disjoint Hamilton cycles.
The following result confirms the  conjecture of Erd\H{o}s.
\begin{theorem}\label{erdos}
Let $T$ be a tournament on $n$ vertices which is chosen uniformly at random. Then a.a.s.~$T$ contains
$\delta^0(T)$ edge-disjoint Hamilton cycles.
\end{theorem}
Note that Theorem~\ref{erdos} is equivalent to the following statement: Consider the complete graph on $n$ vertices and orient each
edge randomly (where the probability for each of the two possible directions is $1/2$), independently of all other edges.
Let $T$ denote the resulting tournament. Then a.a.s.~$T$ contains $\delta^0(T)$ edge-disjoint Hamilton cycles.
To prove Theorem~\ref{erdos}, we find a $\delta^0(T)$-regular oriented spanning subgraph $T'$ of $T$ and apply Theorem~\ref{orientcor} to $T'$
to obtain a Hamilton decomposition of $T'$.

A similar phenomenon occurs also for undirected random graphs.
Let $G_{n,p}$ denote the binomial random graph with edge probability $p$.
Bollob\'as and Frieze~\cite{BF85}, showed that
a.a.s.~$G_{n,p}$ contains $\lfloor \delta(G_{n,p})/2 \rfloor$ edge-disjoint Hamilton cycles in the range of $p$ where the minimum degree $\delta(G_{n,p})$
is a.a.s.~bounded.
Frieze and Krivelevich~\cite{FK08} conjectured that this result extends to the entire range of edge probabilities $p$.
Partial results were proved in several papers, e.g.~\cite{FK08, BKS, AHDoRG, KnoxKO, KrS}. 
In particular, Knox, K\"uhn and Osthus~\cite{KnoxKO} confirmed the conjecture in the range when $(\log n)^{50}/n\le p\le 1-n^{-1/4} (\log n)^9$
and shortly afterwards Krivelevich and Samotij~\cite{KrS} covered the range when $\log n/n\le p\le n^{-1+\eps}$. So in combination with~\cite{BF85},
this implies that the conjecture remains open only in the case when $p$ tends to  $1$ fairly quickly, i.e.~when $p \ge 1-n^{-1/4} (\log n)^9$.
This special case follows without too much work from Theorem~\ref{reggraph} (in a similar way as Theorem~\ref{erdos}).
Altogether, this gives the following result, which confirms the conjecture of Frieze and Krivelevich.
We emphasize that the main contribution to Theorem~\ref{Gnp} comes from the results in~\cite{KnoxKO,KrS}.
\begin{theorem}\label{Gnp}
For any $p=p(n)$, a.a.s.~$G_{n,p}$ contains $\lfloor \delta(G_{n,p})/2\rfloor$ edge-disjoint Hamilton cycles.
\end{theorem}

\subsection{Quasi-random graphs}
Robust expansion is a significant weakening of the extremely well-studied notion of quasi-randomness.
The latter is usually defined in terms of the eigenvalues of a graph. More precisely, 
given a graph $G$, let $\lambda_1\ge \lambda_2\ge \dots\ge \lambda_n$ denote the eigenvalues of its adjacency matrix.
The \emph{second eigenvalue} of $G$ is $\lambda(G):=\max_{i\ge 2} |\lambda_i|$.
We say that a graph $G$ is an \emph{$(n,d,\lambda)$-graph} if it has $n$ vertices, is $d$-regular and $\lambda(G)\le \lambda$. 
It is well known that if $\lambda$ is much smaller than $d$, then such a graph $G$ has strong expansion properties.
This means that we can apply Theorem~\ref{decomp} to quasi-random graphs in order to obtain the following result.
\begin{theorem} \label{corquasi}
For every $\alpha>0$ there exist $\eps>0$ and $n_0\in\mathbb{N}$ such that the following holds.
Let $G$ be any $(n,d,\lambda)$-graph on $n\ge n_0$ vertices with $\lambda\le \eps n$ and such that $d\ge \alpha n$ is even.
Then $G$ has a Hamilton decomposition.
\end{theorem}
Since it is well known that Paley graphs satisfy strong quasi-randomness conditions (see e.g.~\cite{Bollobasbook}), Theorem~\ref{corquasi} 
generalizes (for large $n$) a recent result of 
Alspach, Bryant and Dyer~\cite{Paley} that every Paley graph has a Hamilton decomposition.
Unsurprisingly, dense random regular graphs also satisfy strong quasi-random properties (see~\cite{KSVW} for precise results). So 
Theorem~\ref{corquasi} implies that a.a.s.~random $r$-regular graphs have a Hamilton decomposition if 
$r$ is even and linear in $n$.
The case when $r$ is bounded has received much attention and was settled by Kim and Wormald~\cite{KW01}.
Frieze and Krivelevich~\cite{fk} proved that $\eps$-regular graphs which are almost regular have an approximate Hamilton decomposition.
As graphs which are regular and $\eps$-regular satisfy the conditions of Theorem~\ref{corquasi}, their result can be viewed as an 
approximate version of Theorem~\ref{corquasi}.

\subsection{Partite tournaments}

Jackson~\cite{jackson} posed the following bipartite version of Kelly's conjecture
(both versions are also discussed e.g.~by Bondy~\cite{bondy}).
Here a \emph{bipartite tournament} is an orientation of a complete bipartite graph.
\begin{conj}[Jackson] \label{kellybip}
Every regular bipartite tournament has a Hamilton decomposition.
\end{conj}
This does not seem to follow from our results. However, a $k$-partite version for $k \ge 4$ follows from Theorems~\ref{orientcor} and~\ref{decomp}.
For this, we define a \emph{regular $k$-partite tournament} to be 
an orientation of a complete $k$-partite graph with equal size vertex classes in which the indegree of every vertex equals its outdegree.
\begin{cor} \label{multicor}
For every $k \in \mathbb{N}$ with $k \ge 4$, there is an $n_0 \in \mathbb{N}$ so that every regular $k$-partite tournament $G$ on $n\ge n_0$ vertices has a Hamilton decomposition.
\end{cor}
Note that a regular $k$-partite tournament has degree $(k-1)n/2k$, so for $k \ge 5$ it satisfies the conditions of Theorem~\ref{orientcor}, and so
Corollary~\ref{multicor} follows immediately.
For $k=4$, one cannot apply Theorem~\ref{orientcor}, but one can show  that a regular $4$-partite tournament 
is a robust outexpander (this follows by analyzing e.g.~the proof of Lemma~13.1.~in~\cite{monster}).%
    \COMMENT{See the end of the latex file for the details.} 
Thus one can apply Theorem~\ref{decomp} directly in this case to obtain a Hamilton decomposition.
A regular $3$-partite tournament is not necessarily a robust outexpander (e.g.~it could be a blow-up of an oriented triangle).
But we conjecture that Corollary~\ref{multicor} can be extended to all $k \ge 2$.
\begin{conj}
Every regular tripartite tournament has a Hamilton decomposition.
\end{conj}%Note that it is not even obvious that a bipartite tournament has a single Hamilton cycle --
%this was proved in~\cite{jackson}. 

The corresponding question for undirected complete $k$-partite graphs was settled in the affirmative by 
Auerbach and Laskar~\cite{auerbach}, as well as Hetyei~\cite{hetyei}.
A digraph version was proved by Ng~\cite{Ng}.

\subsection{Algorithmic aspects}
Robust expansion has turned out to be a useful and natural concept -- not just for Hamilton decompositions (see e.g.~\cite{KT}).
This raises the question of whether this property can be recognized efficiently. 
The following result answers this question in the affirmative, as long as one does not need to know the exact expansion parameters.
This is a similar situation as for the concept of $\eps$-regularity. 
Note that for both concepts, the exact parameters are not relevant for any of the applications.
\begin{theorem}\label{algo}
Given $0<\nu\le \tau<1$, there is an algorithm which in time polynomial in~$n$ either decides that a digraph $G$ on $n$ vertices
is not a robust $(\nu,\tau)$-outexpander or decides that $G$ is a robust $(\nu^3/2560,6\tau)$-outexpander.
A similar statement holds for (undirected) robust expanders.
\end{theorem}
The proof of Theorem~\ref{algo} proceeds by showing that a (di-)graph $G$ is a robust (out-)\-expander if and only if the reduced (di-)graph $R$
obtained from an application of Szemer\'edi's regularity lemma is one.
Similar ideas are also used to make the proof of Theorem~\ref{undir_decomp} algorithmic in Section~\ref{sec:algo}.

It would be interesting to know if one could obtain a characterization of robust (out-)expansion without considering $R$,
but perhaps in terms of the eigenvalues of $G$ and/or the codegrees (i.e.~common neighbourhoods) of pairs of vertices.
Such a connection would appear to be quite natural because of the well-known connection between expansion, eigenvalue separation and 
codegree bounds.

\subsection{Organization of the paper}
In Section~\ref{sec:undir}, we derive Theorem~\ref{undir_decomp} from its directed version (Theorem~\ref{decomp}).
We then deduce Theorems~\ref{reggraph} and~\ref{mingraph} from Theorem~\ref{undir_decomp}.
In Section~\ref{sec:dir}, we use Theorem~\ref{decomp} to prove Theorem~\ref{mindigraph}.
We then consider (quasi-)random structures in Section~\ref{sec:tourn} to prove
Theorems~\ref{erdos},~\ref{Gnp} and~\ref{corquasi}.
We give an algorithm for checking robust expansion in Section~\ref{sec:algo}
as well as an algorithmic proof of Theorem~\ref{undir_decomp}.
%Both proofs are based on Szemer\'edi's regularity lemma
%%%%%%%%%%%%%%%%%%%%%%%%%%%%%%%%%%%%%%%%%%%%%%%%%%%%%%%%%%%%%%%%%%%%%%%%%%%%%%%%%%%%%%%%%%%%%%

\section{Notation and a Chernoff type bound}

Given a graph or digraph $G$, we denote its vertex set by~$V(G)$ and its edge set by $E(G)$. We also write $|G|$ for the number of vertices in $G$.

If $G$ is an undirected graph, we write $\delta(G)$ for the minimum degree of $G$ and $\Delta(G)$ for its maximum degree. Given $X,Y\subseteq V(G)$,
we denote the number of all edges of $G$ with one endvertex in $X$ and the other in $Y$ by $e_G(X,Y)$.
An \emph{$r$-factor} of a graph $G$ is a spanning subgraph $H$ of $G$ in which every vertex has degree~$r$.

%If $X\cap Y=\emptyset$ we write
%$G[X,Y]$ for the bipartite subgraph of $G$ whose vertex classes are $X$ and $Y$ and whose edges are all the edges of $G$ between $X$ and $Y$.

If $G$ is a digraph, we write $xy$ for an edge which is directed from $x$ to $y$. The \emph{outneighbourhood $N^+_G(x)$} of a vertex $x$ is the set of all those
vertices $y$ for which $xy\in E(G)$. The \emph{inneighbourhood $N^-_G(x)$} of $x$ is the set of all those
vertices $y$ for which $yx\in E(G)$. We write $d^+_G(x):=|N^+_G(x)|$ for the \emph{outdegree of $x$} and $d^-_G(x):=|N^-_G(x)|$ for the \emph{indegree of $x$}.
We write $\delta^+(G):=\min_{x\in V(G)} d^+_G(x)$ for the \emph{minimum outdegree} of $G$,
$\delta^-(G):=\min_{x\in V(G)} d^-_G(x)$ for its \emph{minimum indegree}, $\delta^0(G):=\min\{\delta^+(G),\delta^-(G)\}$ for its
\emph{minimum semidegree} and $\Delta^0(G):=\max\{\delta^+(G),\delta^-(G)\}$ for the
\emph{maximum semidegree} of $G$. Given $X,Y\subseteq V(G)$,
we denote the number of all edges of $G$ with initial vertex in $X$ and final vertex $Y$ by $e_G(X,Y)$.
%If $X\cap Y=\emptyset$ we write
%$G[X,Y]$ for the bipartite subdigraph of $G$ whose vertex classes are $X$ and $Y$ and whose edges are all the edges of $G$
%with initial vertex in $X$ and final vertex $Y$.
In all these definitions, we omit the subscript $G$ if the graph or digraph $G$ is clear from the context.
An \emph{$r$-factor} of a digraph $G$ is a spanning subdigraph $H$ of $G$ with
$d^+_H(x)=r=d^-_H(x)$ for every $x\in V(G)$. 

The constants in the hierarchies used to state our results have to be chosen from right to left.
More precisely, if we claim that a result holds whenever $0<1/n\ll a\ll b\ll c\le 1$ (where $n$ is the order of the graph or digraph), then this means that
there are non-decreasing functions $f:(0,1]\to (0,1]$, $g:(0,1]\to (0,1]$ and $h:(0,1]\to (0,1]$ such that the result holds
for all $0<a,b,c\le 1$ and all $n\in \mathbb{N}$ with $b\le f(c)$, $a\le g(b)$ and $1/n\le h(a)$. 
We will not calculate these functions explicitly.

We will often use the following Chernoff type bound (see e.g.~Theorem~2.1 and Corollary 2.2 in~\cite{JLR}).

\begin{prop}\label{chernoff}
Suppose $X$ has binomial distribution and $0<a<1$. Then
$$
\mathbb{P}(X \le (1-a)\mathbb{E}X) \le e^{-\frac{a^2}{3}\mathbb{E}X} \ \ \ \text{and} \ \ \
\mathbb{P}(X \ge (1+a)\mathbb{E}X) \le e^{-\frac{a^2}{3}\mathbb{E}X}.
$$
\end{prop}

%%%%%%%%%%%%%%%%%%%%%%%%%%%%%%%%%%%%%%%%%%%%%%%%%%%%%%%%%%%%%%%%%%%%%%%%%%%%%%%%

\section{Edge-disjoint Hamilton cycles in undirected graphs} \label{sec:undir}
The purpose of this section is to prove Theorem~\ref{undir_decomp}. We will then use it to derive
Theorems~\ref{reggraph} and~\ref{mingraph}.

\subsection{Regular orientations and Theorem~\ref{undir_decomp}}
To prove Theorem~\ref{undir_decomp}, we will show that the edges of every $r$-regular (undirected) robust expander can
be oriented such that the oriented graph $G^{\rm orient}$ thus obtained from $G$ is a $r/2$-regular robust outexpander
(see Lemma~\ref{reggraphorient}).
Theorem~\ref{decomp} then implies that $G^{\rm orient}$ has a Hamilton decomposition, which clearly corresponds to a
Hamilton decomposition of $G$. 

In order to prove Lemma~\ref{reggraphorient}, we first show that there is some orientation
$G'$ of $G$ such that $G'$ is a robust outexpander (see Lemma~\ref{orientexp}). It is not hard to check that a random orientation
will satisfy this property. However, it will only be almost regular, and not regular. We then show that every robust
outexpander (and thus also $G'$) contains a sparse \emph{regular} spanning subdigraph $G^*$ which is still a robust outexpander
(see Lemma~\ref{regexpslice}). It now remains
to find a regular orientation $G^\diamond$ of $G\setminus E(G^*)$. Then $G^{\rm orient}:=G^*\cup G^\diamond$ is a $r/2$-regular
orientation of $G$ and, since $G^*$ is
a robust outexpander, $G^{\rm orient}$ is still a robust outexpander.

The proofs of Lemmas~\ref{orientexp} and~\ref{regexpslice} (and thus of Lemma~\ref{reggraphorient}) are not algorithmic.
In Section~\ref{sec:orient} we will also prove an algorithmic version of Lemma~\ref{reggraphorient}
(which leads to an algorithmic proof of Theorem~\ref{undir_decomp}) based on Szemer\'edi's regularity lemma.
However, this does not mean that the corresponding part of the current subsection can be omitted:
Lemma~\ref{regexpslice} (and thus also Lemmas~\ref{slicerobust} and~\ref{regrobust}) will be used again in our proof of Theorem~\ref{mindigraph}.
Moreover, Lemmas~\ref{orientexp} and~\ref{regexpslice} are also crucial ingredients in~\cite{KLOmindeg}.
So it is only the very short derivation of Lemma~\ref{reggraphorient} itself that might be considered redundant.
We have included it for completeness, as the underlying observation that addition of edges preserves robust expansion  
seems to be very useful (as is also illustrated in~\cite{KLOmindeg, FKS}).

As described above,  we first show that every (undirected) robust expander has an almost regular orientation which is
a robust outexpander. 

\begin{lemma}\label{orientexp}
Suppose that $0<1/n\ll \eta\ll \nu,\tau,\alpha<1$.
Suppose that $G$ is a robust $(\nu,\tau)$-expander on $n$ vertices with $\delta(G)\ge \alpha n$.
Then one can orient the edges of $G$ in such a way that the oriented graph $G'$ thus obtained from $G$ satisfies the following properties:
\begin{itemize}
\item[{\rm (i)}]  $G'$ is a robust $(\nu/4,\tau)$-outexpander.
\item[{\rm (ii)}] $d^+_{G'}(x)=(1\pm \eta)d_G(x)/2$ and $d^-_{G'}(x)=(1\pm \eta)d_G(x)/2$ for every vertex $x$ of $G$.
\end{itemize}
\end{lemma}
\proof
Orient each edge $xy$ of $G$ randomly (where the probability for each of the two possible directions  is $1/2$), independently of all other edges.
Let $G'$ denote the oriented graph obtained in this way. Using Proposition~\ref{chernoff} it is easy to show that (ii) fails with probability at most $1/4$.

So consider any $S\subseteq V(G)$ with $\tau n\le |S|\le (1-\tau)n$.
Let $U':=RN^+_{\nu/4,G'}(S)$ and $U:=RN_{\nu,G}(S)$.
Let $U^*$ be any subset of $U$ of size $\nu n/4$.
Consider any $u\in U^*$ and let $X:=N^-_{G'}(u)\cap (S\setminus U^*)$. Note that
$$\ex X\ge \frac{\nu n-|U^*|}{2}=\frac{3\nu n}{8}.$$
So Proposition~\ref{chernoff} implies that
\begin{align*}
\pr(X\le \nu n/4)& \le \pr(X \le 2\ex X/3) \le e^{-\ex X/27}\le e^{-\nu n/100}.
\end{align*}
Thus the probability that all vertices in $U^*$ have at most $\nu n/4$ inneighbours in $S\setminus U^*$
(in the oriented graph $G'$)
is at most $e^{-\nu^2 n^2/400}$ (note these events are independent for different vertices in $U^*$ as we are only considering the
inneighbours in~$S\setminus U^*$ rather than all inneighbours in~$S$).
But if $|U'|\le |S|+\nu n/4$, then there exists a set $U^*\subseteq U$ of size $\nu n/4$ such that every vertex in $U^*$
has at most $\nu n/4$ inneighbours in $S\setminus U^*$. So the probability that $U'$ has size at most $|S|+\nu n/4$ is at most
$$
\binom{|U|}{\nu n/4} e^{-\nu^2 n^2/400} \le 2^{n} e^{-\nu^2 n^2/400} \le e^{-\nu^2 n^2/401}.
$$
Summing over all possible sets $S$ shows that the probability that (i) fails is also at most $1/4$.
\endproof

The next lemma shows that the edges of any robust outexpander can be split  in such a way that
the two digraphs thus obtained are still robust outexpanders.

\begin{lemma} \label{slicerobust}
Suppose that $0<1/n\ll \eta \ll \nu,\tau,\alpha,\lambda, 1-\lambda<1$.
Let~$G$ be a digraph on~$n$ vertices with $\delta^0(G)\ge \alpha n$ which is a robust $(\nu,\tau)$-outexpander.
Then $G$ can be split into two edge-disjoint spanning subdigraphs $G_1$ and $G_2$ such that the following two properties hold.
\begin{itemize}
\item[(i)] $d^+_{G_1}(x)=(1\pm \eta)\lambda d^+_G(x)$ and $d^-_{G_1}(x)=(1\pm \eta)\lambda d^-_G(x)$ for every $x\in V(G)$. 
\item[(ii)] $G_1$ is a robust $(\lambda\nu/2,\tau)$-outexpander and $G_2$ is a robust $((1-\lambda)\nu/2,\tau)$-outexpander.
\end{itemize} 
\end{lemma}
\proof
Consider a random partition of the edges of $G$ into two subdigraphs $G_1$ and $G_2$ where an edge is included into $G_1$
with probability $\lambda$ and into $G_2$ with probability $1-\lambda$ (independently of all other edges).
%Consider a random partition of the edges of $G$ into $s$ parts, where for each edge the probability of adding $G_j$
%is $1/s$ (independently of all other edges).
Proposition~\ref{chernoff} immediately implies that the probability that (i) fails is at most $1/4$.

The argument for (ii) is similar to the proof of Lemma~\ref{orientexp}. Consider any $S\subseteq V(G)$ with $\tau n\le |S|\le (1-\tau)n$.
Let $U:=RN^+_{\nu,G}(S)$, $U_1:=RN^+_{\lambda\nu /2,G_1}(S)$ and $U_2:=RN^+_{(1-\lambda)\nu /2,G_2}(S)$.
Consider any $u\in U$. The expected number of edges from $S$ to $u$ in $G_1$ is at least $\nu \lambda n$.
So Proposition~\ref{chernoff} implies that the probability that $u$ does not lie in $U_1$ is at most $e^{-\nu \lambda n/12}$.
Let $U'$ be any subset of $U$ of size $\nu n/2$.
Then the probability that no vertex of $U'$ lies in $U_1$ is at most $e^{-\nu^2\lambda n^2/24}$
(since these events are independent for different vertices of $U'$).
So the probability that $U_1$ has size at most $|S|+\nu n/2$ is at most
$$
\binom{|U|}{\nu n/2} e^{-\nu^2\lambda n^2/24} \le 2^{n} e^{-\nu^2\lambda n^2/24} \le e^{-\nu^2\lambda n^2/25}.
$$
One can use a similar argument to show that the probability that $U_2$ has size at most $|S|+\nu n/2$ is at most
$e^{-\nu^2(1-\lambda) n^2/25}$. Summing over all possible sets $S$ shows that the probability that (ii) fails is also at most $1/4$.
\endproof

The next lemma (which was proved in~\cite{monster}) will be used in the proof of Lemma~\ref{regexpslice} to turn
a sparse digraph $G_1$ into a regular one by adding a digraph $G'_2$ whose degree sequence complements that of $G_1$.
The lemma guarantees the existence of such a digraph $G'_2$ in any robust outexpander, provided that the degrees of $G'_2$ are
within a certain range.

\begin{lemma}\label{regrobust}
Suppose that $0<1/n\ll\eps,\xi\ll \nu \le \tau \ll \alpha<1$. Let~$G$ be a digraph on~$n$ vertices with
$\delta^0(G)\ge \alpha n$ which is a robust $(\nu,\tau)$-outexpander. For every vertex $x$ of $G$ let $n^+_x,n^-_x\in\mathbb{N}$
be such that $(1-\eps)\xi n\le n^+_x, n^-_x\le (1+\eps)\xi n$ and such that $\sum_{x\in V(G)} n^+_x=\sum_{x\in V(G)} n^-_x$.
Then $G$ contains a spanning subdigraph $G'$ such that $d^+_{G'}(x)=n^+_x$ and $d^-_{G'}(x)=n^-_x$ for every $x\in V(G)$.
\end{lemma}

In order to prove Lemma~\ref{regexpslice}, we first apply Lemma~\ref{slicerobust} to split the edges of $G$ into two digraphs $G_1$ and $G_2$ which
are both still robust outexpanders, where $G_1$ is very sparse. We then use Lemma~\ref{regrobust}
to find a subdigraph $G'_2$ inside $G_2$ whose degree sequence complements that
of $G_1$. So $G_1\cup G'_2$ will then be a \emph{regular} robust outexpander.

\begin{lemma}\label{regexpslice}
Suppose that $0<1/n\ll \nu'\ll \xi\ll\nu\le \tau\ll \alpha<1$. Let $G$ be a robust $(\nu,\tau)$-outexpander on $n$ vertices
with $\delta^0(G)\ge \alpha n$. Then $G$ contains a $\xi n$-factor which is still a robust $(\nu',\tau)$-outexpander.
\end{lemma}
\proof
Choose new constants $\eta, \lambda$ such that $0<1/n\ll \eta,\nu' \ll \lambda \ll \xi$.
Apply Lemma~\ref{slicerobust} to partition $G$ into edge-disjoint
spanning subdigraphs $G_1$ and $G_2$ such that the following properties are satisfied:
\begin{itemize}
\item $d^+ _{G_1}(x)=(1\pm \eta)\lambda d^+_G(x)$ and $d^-_{G_1}(x)=(1\pm \eta)\lambda d^-_G(x)$ for every $x\in V(G)$.
\item $G_1$ is a robust $(\lambda\nu /2,\tau)$-outexpander and $G_2$ is a robust $((1-\lambda)\nu /2,\tau)$-outexpander.
\end{itemize}
For every $x\in V(G)$, let $n^\pm_x:=\xi n-d^\pm_{G_1}(x)$.
Note that
$$ (1-\sqrt{\lambda})\xi n\le \xi n- 2\lambda n\le n^+_x,n^-_x\le \xi n.
$$
Moreover, $\delta^0(G_2)\ge \alpha n-(1+ \eta)\lambda n\ge \alpha n/2$.
Thus we can apply Lemma~\ref{regrobust} to $G_2$ (with $\sqrt{\lambda}$, $(1-\lambda)\nu /2$, $\alpha/2$ playing the roles of $\eps$, $\nu$, $\alpha$)
to obtain a spanning subdigraph $G'_2$ of $G_2$ such that $d^\pm_{G'_2}(x)=n^\pm_x$ for every $x\in V(G)$. Then $G^*:=G_1\cup G'_2$
is a $\xi n$-factor of $G$. Moreover, $G^*$ is a robust $(\nu',\tau)$-outexpander since it contains the robust $(\lambda\nu /2,\tau)$-outexpander
$G_1$ (and since $\nu'\le \lambda\nu /2$).
\endproof

We also use the following classical result of Petersen.

\begin{theorem}\label{petersen}
Every regular graph of positive even degree contains a $2$-factor.
\end{theorem}

We can now combine Lemmas~\ref{orientexp} and~\ref{regexpslice} as well as Theorem~\ref{petersen} to prove the main ingredient
for our proof of Theorem~\ref{undir_decomp}.

\begin{lemma}\label{reggraphorient}
Suppose that $0<1/n\ll \nu'\ll \nu\le \tau\ll \alpha<1$.
Let $G$ be an $r$-regular graph on $n$ vertices such that $r\ge \alpha n$ is even and $G$ is a robust $(\nu,\tau)$-expander.
Then one can orient the edges of $G$ in such a way that the oriented graph $G^{\rm orient}$ thus obtained from $G$ is an
$r/2$-regular robust $(\nu',\tau)$-outexpander.
\end{lemma}
\proof
Choose a new constant $\xi$ such that $\nu' \ll \xi \ll \nu$.
Apply Lemma~\ref{orientexp} to find an orientation $G'$ of $G$ such that $G'$ is a robust $(\nu/4,\tau)$-outexpander
and such that $\delta^0(G')\ge \alpha n/3$ (say).
Now apply Lemma~\ref{regexpslice} to find a $\xi n$-factor $G^*$ of $G'$ which is still a robust $(\nu',\tau)$-outexpander.
Let $H$ be the undirected graph obtained from $G$ by deleting all the edges in $G^*$.
Then $H$ is $(r-2\xi n)$-regular. Thus Petersen's theorem (Theorem~\ref{petersen}) implies that $H$ has a decomposition into edge-disjoint 2-factors.
Orient each cycle in these 2-factors consistently. This gives an orientation of $H$ which is $(r-2\xi n)/2$-regular.
Adding the edges of $G^*$ gives an orientation $G^{\rm orient}$ of $G$ which is $r/2$-regular. Moreover, $G^{\rm orient}$
is a robust $(\nu',\tau)$-outexpander since it contains the robust $(\nu',\tau)$-outexpander $G^*$.
\endproof

We can now prove an undirected analogue of Theorem~\ref{decomp} on Hamilton decompositions of robust expanders.

\removelastskip\penalty55\medskip\noindent{\bf Proof of Theorem~\ref{undir_decomp}. }
Let $\tau^*:=\tau(\alpha)$, where $\tau(\alpha)$ is as defined in Theorem~\ref{decomp}.
Choose a new constant $\tau$ such that $0<  \tau\ll \alpha,\tau^*$.
Note that whenever $\nu'\le \nu$, every robust $(\nu,\tau)$-outexpander is also a robust $(\nu',\tau)$-outexpander.
So we may assume that $0 \ll \nu \ll \tau$.
Now choose $n_0 \in \mathbb{N}$ and $\nu'$ such that $1/n_0 \ll \nu'\ll \nu$.
($\tau$ and $n_0$ will be the constants returned by Theorem~\ref{undir_decomp}.)
Now let $G$ be a graph satisfying the conditions of Theorem~\ref{undir_decomp} with these parameters, 
i.e.~$G$ is an $r$-regular robust $(\nu,\tau)$-expander on $n\ge n_0$ vertices such that $r\ge \alpha n$ is even.
Apply Lemma~\ref{reggraphorient} to obtain an $r/2$-regular orientation $G^{\rm orient}$ of $G$ such that $G^{\rm orient}$ is a
robust $(\nu',\tau)$-outexpander (and thus also a robust $(\nu',\tau^*)$-outexpander).
Now apply Theorem~\ref{decomp}  to find a Hamilton decomposition of $G^{\rm orient}$
(with $\nu'$ playing the role of $\nu$ in Theorem~\ref{decomp}).
Clearly, this Hamilton decomposition corresponds to a Hamilton decomposition of $G$.
\endproof

\subsection{Theorems~\ref{reggraph} and~\ref{mingraph}}
Our next aim is to derive Theorem~\ref{reggraph} from Theorem~\ref{undir_decomp}.
For this, it suffices to show that every graph of minimum degree a little larger than $n/2$ is a robust expander.
%\begin{lemma}\label{reggraphexpander}
%Suppose that $0<\nu\le \tau\le \eps<1$ are such that $\eps\ge 2\nu/\tau$.
%Let $G$ be a graph on $n$ vertices with minimum degree $\delta(G)\ge (1/2+\eps)n$.
%Then $G$ is a robust $(\nu,\tau)$-expander.
%\end{lemma}
%\proof
%Consider any set $S\subseteq V(G)$ with $\tau n\le |S|\le  (1-\tau)n$.
%Let $RN:=RN_{\nu,G}(S)$. We have to show that $|RN|\ge |S|+\nu n$. Suppose first that $|S|\ge n/2$.
%Then every vertex of $G$ has at least $\eps n\ge \nu n$ inneighbours in $S$. So $RN=V(G)$. So we may assume that
%$|S|\le n/2$. 
%For any sets of vertices $X,Y$, let $e'(X,Y):=\sum_{x \in X} |N(x) \cap Y|$. But
%\begin{align*}
%(1/2+\eps)n|S| \le e'(S,V(G)) & =e'(S,RN)+e'(S,V(G)\setminus RN)\le |S||RN|+\nu n^2\\
%& \le |S||RN|+\frac{\nu}{\tau}n|S|.
%\end{align*}
%So $|RN|\ge (1/2+\eps-\nu/\tau)n\ge (1+\eps)n/2\ge |S|+\nu n$, as required.
%\endproof 
To prove this, we will use the following analogue of this result for regular digraphs, which was proved as Lemma~13.2 in~\cite{monster}.
Its proof follows easily from the definition of robust outexpansion.
We will also use Lemma~\ref{regdiexpander} in the proof of Theorem~\ref{mindigraph}.

\begin{lemma}\label{regdiexpander}
Suppose that $0<\nu\le \tau\le \eps<1$ are such that $\eps\ge 2\nu/\tau$.
Let $G$ be a digraph on $n$ vertices with minimum semidegree $\delta^0(G)\ge (1/2+\eps)n$.
Then $G$ is a robust $(\nu,\tau)$-outexpander.
\end{lemma}

\begin{lemma}\label{reggraphexpander}
Suppose that $0<\nu\le \tau\le \eps<1$ are such that $\eps\ge 2\nu/\tau$.
Let $G$ be a graph on $n$ vertices with minimum degree $\delta(G)\ge (1/2+\eps)n$.
Then $G$ is a robust $(\nu,\tau)$-expander.
\end{lemma}
\proof
Let $G'$ be the digraph obtained from $G$ by replacing every (undirected) edge $xy$ of $G$ by two
directed edges $xy$ and $yx$. Then $\delta^0(G')\ge (1/2+\eps)n$ and so Lemma~\ref{regdiexpander} implies
that $G'$ is a robust $(\nu,\tau)$-outexpander. But this implies that $G$ is a robust $(\nu,\tau)$-expander.
\endproof

Theorem~\ref{reggraph} (which guarantees a Hamilton decomposition of dense regular graphs) 
is now an immediate consequence of Lemma~\ref{reggraphexpander} and Theorem~\ref{undir_decomp}.

\removelastskip\penalty55\medskip\noindent{\bf Proof of Theorem~\ref{reggraph}. }
Let $\tau^*:=\tau(1/2)$, where $\tau(1/2)$ is as defined in Theorem~\ref{undir_decomp}.
Choose new constants $n_0\in\mathbb{N}$ and $\nu, \tau$ such that $0<1/n_0\ll \nu\ll \tau\le \eps,\tau^*$. 
Now let $G$ be a graph satisfying the conditions of Theorem~\ref{reggraph}, i.e.~$G$ is an $r$-regular graph on $n \ge n_0$ vertices, 
where $r \ge (1/2+\eps) n$ is even.
Then Lemma~\ref{reggraphexpander} implies that $G$ is a robust $(\nu,\tau)$-expander (and so also a
robust $(\nu,\tau^*)$-expander).
Thus Theorem~\ref{undir_decomp} implies that $G$ has a Hamilton decomposition.
\endproof

Finally, we can use Theorem~\ref{reggraph} to derive Theorem~\ref{mingraph}
(which concerns edge-disjoint Hamilton cycles in graphs of given minimum degree).

\removelastskip\penalty55\medskip\noindent{\bf Proof of Theorem~\ref{mingraph}. }
Choose $n_0\in\mathbb{N}$ and an additional constant $\eps'$ such that $1/n_0 \ll \eps' \ll\eps$.
Consider any graph $G$ on $n\ge n_0$ vertices with minimum degree $\delta\ge (2-\sqrt{2}+\eps)n$. Then Theorem~\ref{CKOthm} implies that
$${\rm reg}_{\rm even}(G)\ge {\rm reg}_{\rm even}(n,\delta)\ge g_{\rm even}(n,\delta).$$
(Recall that $g_{\rm even}(n,\delta)$ is the largest even integer $r$ with  $r \le (\delta+\sqrt{n(2\delta-n)})/2$.) 
Note that the only solution to  $a+ \sqrt{2a-1}=1$%
\COMMENT{$a+ \sqrt{2a-1}=1$ iff $2a-1 =(1-a)^2=1-2a+a^2$ iff $a^2-4a = -2$ iff $(a-2)^2=-2+4=2$}
 with $0 \le a \le 1$ is $a=2 - \sqrt{2}$.
So the lower bound on $\delta$ implies that $g_{\rm even}(n,\delta)\ge (1/2+\eps')n$.
Let $G'$ be a regular subgraph of $G$ of degree ${\rm reg}_{\rm even}(G)$. Since
${\rm reg}_{\rm even}(G)\ge g_{\rm even}(n,\delta)\ge (1/2+\eps')n$, Theorem~\ref{reggraph} implies that
$G'$ has a decomposition into ${\rm reg}_{\rm even}(G)/2$ edge-disjoint Hamilton cycles.
\endproof
%%%%%%%%%%%%%%%%%%%%%%%%%%%%%%%%%%%%%%%%%%%%%%%%%%%%%%%%%%%%%%%%%%%%%%%%%%

\section{Edge-disjoint Hamilton cycles in digraphs of large minimum semidegree} \label{sec:dir}

The purpose of this section is to prove Theorem~\ref{mindigraph}. To do this, we will first determine ${\rm reg}^{\rm dir}(n,\delta)$
for all $n/2\le \delta< n$,
i.e.~the largest number $r$ such that every digraph on $n$ vertices of minimum semidegree $\delta$ contains an $r$-factor.
Recall that ${\rm reg}(G)$ denotes the largest degree of a regular spanning subdigraph of $G$.

\begin{lemma}\label{factor} Suppose that $\delta,n\in\mathbb{N}$ are such that  $n/2\le \delta<n$. Let
$$
 r^*:= \frac{\delta+\sqrt{n(2\delta-n)+\1}}{2},
\ \ \ \text{where}  \ \ \
{\bf 1}_{n\not\equiv \delta}:=\begin{cases} 0 & \mbox{if } n\equiv \delta \mod 2\\
1 & \mbox{if } n\not \equiv \delta \mod 2.\end{cases}
$$
\begin{itemize}
\item[\rm{(i)}] Let $G$ be a digraph on $n$ vertices with minimum semidegree $\delta$. Then
${\rm reg}(G) \ge \lfloor r^*\rfloor$.
\item[\rm{(ii)}] There is a digraph $G$ on $n$ vertices with minimum semidegree $\delta$
such that ${\rm reg}(G) = \lfloor r^*\rfloor$. 
\end{itemize}
In particular, ${\rm reg}^{\rm dir}(n,\delta)=\lfloor r^*\rfloor$. 
\end{lemma}
\proof
Let $r:=\lfloor r^* \rfloor$. In order to prove~(i), we have to show that $G$ contains
an $r$-factor.%
   \COMMENT{Since every regular digraph contains a $1$-factor, it follows that $G$ will also contain an $r'$-factor for all $r'<r$.}
To do this, our aim is to apply the Max-Flow-Min-Cut theorem.
So let $H$ be the (unoriented) bipartite graph whose vertex classes $A$ and $B$ are both copies of $V(G)$ and in which
$a\in A$ is joined to $b\in B$ if $ab$ is a (directed) edge of $G$. Give every edge of $H$ capacity~1. Add a source $u^*$ which is
joined to every vertex $a\in A$ with an edge of capacity $r$. Add a sink $w^*$ which is
joined to every vertex $b\in B$ with an edge of capacity $r$. 
Note that an integer-valued $r n$-flow corresponds to the desired $r$-factor of $G$. Thus by the
Max-Flow-Min-Cut theorem it suffices to show that every cut has capacity at least $r n$.

So consider a minimal cut $\mathcal{C}$. Let $U$ be the set of all those vertices $a\in A$ for which $u^*a\notin \mathcal{C}$.
Similarly, let $W$ be the set of all those vertices $b\in B$ for which $bw^*\notin \mathcal{C}$.
We only consider the case when $|W|\ge |U|$ (the other case is similar).
Let $U':=A\setminus U$ and $W':=B\setminus W$.
Thus the capacity of $\mathcal{C}$ is
$$ c:=r |U'|+e_B(U,W)+ r |W'|.
$$
But $$e_B(U,W)=e_G(U,W)\ge \delta |U|-e_G(U,W')\ge  \delta |U|-|U||W'|=r|U|+|U|(\delta-r-|W'|).$$
Thus $$c\ge r |U'|+r|U|+|U|(\delta-r-|W'|)+ r |W'|=rn+|U|(\delta-r-|W'|)+ r |W'|.$$
So it suffices to show that
\begin{equation}\label{eqUW}
0\le |U|(\delta-r-|W'|)+ r |W'|=|W'|(r-|U|)+|U|(\delta-r).
\end{equation}
This holds if $|U|\le r$. So we may assume that $|U|>r$. Since $|W|\ge |U|$ we have that $|W'|\le |U'|$.
So (\ref{eqUW}) holds if
\begin{align}
0 & \le |U'|(r-|U|)+|U|(\delta-r)=(n-|U|)(r-|U|)+|U|(\delta-r) \nonumber\\
& = nr+|U|^2-|U|(n-\delta+2r)\nonumber\\
& =\left( |U|-\frac{n-\delta+2r}{2}\right)^2-\frac{(n-\delta+2r)^2}{4}+nr\label{eq:flow}.
\end{align}
But writing $$f(x):=nx- \frac{(n-\delta+2x)^2}{4},$$ we have that%
     \COMMENT{\begin{align*}
4f(r^*)& =2n(\delta+\sqrt{n(2\delta-n)+\1})- (n+\sqrt{n(2\delta-n)+\1})^2\\ & =
2n(\delta+\sqrt{n(2\delta-n)+\1})-(n^2+2n\sqrt{n(2\delta-n)+\1}+n(2\delta-n)+\1)=-\1.
\end{align*}}
%$$f(r^*)=\begin{cases} 0 & \mbox{if } n\equiv \delta \mod 2\\
%-1/4 & \mbox{if } n\not \equiv \delta \mod 2\end{cases}.$$
$f(r^*)=-\1 /4$.
Moreover, $f(x)$ is decreasing for all $x\ge \delta/2$ and $\delta/2\le r\le r^*$.
(To check the latter, it is easiest to consider the cases when $\delta=n/2$ and when $\delta>n/2$ separately.)%
     \COMMENT{If $\delta>n/2$ then $r^*\ge (\delta+\sqrt{n})/2$ and so $r\ge \delta/2$.
If $\delta=n/2$ then $r^*=(\delta+\1)/2$ and $n=2\delta$ is even. If $\delta$ is even then
$r^*=\delta/2=r$. If $\delta$ is odd then $r^*=(\delta+1)/2=r$.}
Altogether this shows that
%$$f(r) \ge \begin{cases} 0 & \mbox{if } n\equiv \delta \mod 2\\
%-1/4 & \mbox{if } n\not \equiv \delta \mod 2\end{cases}.$$
$f(r) \ge -\1 /4$. 

On the other hand, $r \in \mathbb{N}$ implies that 
$$\left( |U|-\frac{n-\delta+2r}{2}\right)^2 \ge\begin{cases} 0 & \mbox{if } n\equiv \delta \mod 2\\
1/4 & \mbox{if } n\not \equiv \delta \mod 2\end{cases}.$$
Altogether this shows that (\ref{eq:flow}) holds. This proves~(i).

The proof of~(ii) is similar to the proof of its analogue for undirected graphs (see Section~2 in~\cite{CKO}).
If $\delta=n-1$, then $r^*=n-1$ and so we can take $G$ to be the complete digraph on $n$ vertices.% 
    \COMMENT{In this case $r^*=\frac{n-1+\sqrt{n(n-2)+1}}{2}=\frac{n-1+\sqrt{(n-1)^2}}{2}=n-1$.}
Thus we may assume that $\delta\le n-2$.
Let $$\Delta:=\left\lceil\frac{n+\sqrt{n(2\delta-n)+\1}}{2}\right\rceil.$$ It is easy to check that%
     \COMMENT{Indeed, the first inequality holds if $(2\delta-n)^2\le n(2\delta-n)+\1$.
But $(2\delta-n)^2\le n(2\delta-n)$ is equivalent to $2\delta-n \le n$, which holds. Using that $\delta\le n-2$, we see that the 2nd inequality holds since
$\frac{n+\sqrt{n(2\delta-n)+1}}{2}\le \frac{n+\sqrt{n(n-4)+1}}{2}< \frac{n+\sqrt{(n-2)^2}}{2}=n-1.$}
$\delta < \Delta\le n-1$. (Here we use that $(2\delta-n)^2 < n(2\delta-n)$ for the lower bound and $\delta\le n-2$ for the upper bound.)
Let $A$ be an empty digraph on $n-\Delta$ vertices and let $B$ be
a $(\delta+\Delta-n)$-regular digraph on $\Delta$ vertices.%
    \COMMENT{In~\cite{CKO} we require $\Delta$ to be even to ensure the existence of $B$, but we don't need this for directed graphs.}
Let $G$ be obtained from the disjoint union of $A$ and $B$
by adding all edges from $A$ to $B$ and all edges from $B$ to $A$. Thus $\delta^0(G)=\delta$ and $\Delta^0(G)=\Delta$.
We claim that ${\rm reg}(G)\le \frac{\Delta(\delta+\Delta-n)}{2\Delta-n}$. Indeed, given any $r$-factor $G'$ of $G$,
we have $$
r\Delta=\sum_{x\in B} d^-_{G'}(x)\le \Delta(\delta+\Delta-n)+e_{G'}(A,B)=\Delta(\delta+\Delta-n)+r(n-\Delta),
$$
which implies the claim. But 
$$
\frac{\Delta (\delta + \Delta - n)}{2 \Delta - n} =
\frac{\delta}{2} + \frac{n \delta/2}{2\Delta - n} -
\frac{\Delta(n-\Delta)}{2\Delta - n}.
$$
Let $\eta$ be such that $\Delta=\frac{n+\sqrt{n(2\delta-n)+\1}}{2}+\eta$. So $0\le \eta<1$.
Then
\begin{align*}
\Delta(n - \Delta) &= \left(\frac{n}{2} + \left(\frac{\sqrt{n(2\delta - n)+\1}}{2} + \eta\right) \right)
\left(\frac{n}{2} - \left(\frac{\sqrt{n(2\delta - n)+\1}}{2} + \eta\right) \right)\\
& = \frac{n^2}{4} - \frac{n(2\delta - n)+\1}{4} - \eta\sqrt{n(2\delta - n)+\1} - \eta^2\\
& = \frac{n^2 - n\delta}{2} - \frac{\1}{4}-\eta\sqrt{n(2\delta - n)+\1} - \eta^2.
\end{align*}
Thus
$$
\frac{\Delta (\delta + \Delta - n)}{2 \Delta - n}=\frac{\delta}{2} + \frac{n(2\delta - n) + 2\eta\sqrt{n(2\delta - n)+\1}}{2(2\Delta - n)}
+\frac{\1}{4(2\Delta - n)}+ \frac{\eta^2}{2\Delta - n}.
$$
Since also
$
(2\Delta - n)\sqrt{n(2\delta-n)+\1} = n(2\delta-n)+\1 +
2\eta\sqrt{n(2\delta - n)+\1},
$
we deduce that
\begin{align*}
{\rm reg}(G)& \le \frac{\Delta (\delta + \Delta - n)}{2 \Delta - n}
= \frac{\delta + \sqrt{n(2\delta-n)+\1}}{2}-\frac{\1}{4(2\Delta - n)} +
\frac{\eta^2}{2\Delta - n}\nonumber \\
& =  r^* +\frac{\eta^2 -\1 /4}{2\Delta-n}.            
\end{align*}
If $\eta^2 -\1 /4 \le 0$, this implies that ${\rm reg}(G) \le \lfloor r^* \rfloor$, as required.
So we may assume that  $\eta^2 -\1 /4 >0$. Now recall that $\Delta> \delta$ and hence $2\Delta-n \ge 2(\delta+1) -n \ge 2$.
This means that we may assume that
\begin{equation} \label{eq:fac}
{\rm reg}(G) \le r^* + \eta^2/2- \1/8.
\end{equation}
Let $\eta_*$ be such that $\lceil r^*\rceil=r^*+\eta_*$. So $0\le \eta_*<1$
and it is easy to see that 
$$
\eta = \begin{cases} \eta_*& \mbox{if } n\equiv \delta \mod 2\\
%\eta_*+1/2 \ \ \text{ or } \ \ \eta_*-1/2 & \mbox{if } n\not \equiv \delta \mod 2\end{cases}.
\eta_* \pm 1/2  & \mbox{if } n\not \equiv \delta \mod 2\end{cases}.
$$
Thus if $n\equiv \delta \mod 2$ then $\eta^2/2-\1/8=\eta^2/2=\eta^2_*/2<\eta_*$.
If $n\not \equiv \delta \mod 2$  then
$
\eta^2/2-\1/8 \le \eta^2_*/2+\eta_*/2<\eta_*$. 
This shows that in both cases 
the right hand side of (\ref{eq:fac}) is strictly less than $r^*+\eta_*=\lceil r^*\rceil$,
which in turn implies that ${\rm reg}(G)\le \lfloor r^*\rfloor$. Together with~(i) this shows that
 ${\rm reg}(G)= \lfloor r^*\rfloor$ in all cases, as required.
\endproof

We are now in a position to find (nearly) optimal packings of edge-disjoint Hamilton cycles in digraphs of given minimum semidegree~$\delta$.

\removelastskip\penalty55\medskip\noindent{\bf Proof of Theorem~\ref{mindigraph}. }
We first prove~(i). Lemma~\ref{factor} implies that ${\rm reg}^{\rm dir}(n,\delta)=f(n,\delta)$
whenever $n/2\le\delta< n$. Now suppose that $\delta\ge (1/2+\eps)n$, where $n$ is sufficiently large compared with $1/\eps$.
To show that ${\rm ham}^{\rm dir}(n,\delta)\ge f(n,\delta)-\eps n$, we let
$\tau^*:=\tau(1/5)$, where $\tau(1/5)$ is as defined in Theorem~\ref{decomp}.
Choose new constants $\nu',\nu, \tau,\xi$ such that
$$1/n_0\ll \nu'\ll \xi\ll \nu \ll \tau\le \eps,\tau^*.$$
Consider any digraph $G$ on $n$ vertices with minimum semidegree $\delta$. Then
Lemma~\ref{regdiexpander} implies that $G$ is a robust $(\nu,\tau)$-outexpander.
Apply Lemma~\ref{regexpslice} to obtain a $\xi n$-factor $G_1$ of $G$ such that $G_1$ is still
a robust $(\nu',\tau)$-outexpander (and thus also a robust $(\nu',\tau^*)$-outexpander).
Let $G_2:=G\setminus E(G_1)$ and let $G'_2$ be a spanning subdigraph of $G_2$ which is ${\rm reg}(G_2)$-regular.
Since $\delta^0(G_2)=\delta^0(G)-\xi n=\delta-\xi n$, Lemma~\ref{factor}(i) implies that
${\rm reg}(G_2)\ge f(n,\delta)-\eps n$. Thus $G_1\cup G'_2$ is regular of degree $\xi n+{\rm reg}(G_2)\ge f(n,\delta)-\eps n\ge n/5$
(say). Moreover, $G_1\cup G'_2$ is still a robust $(\nu',\tau^*)$-outexpander
(since it contains the robust $(\nu',\tau^*)$-outexpander $G_1$).
So we may apply Theorem~\ref{decomp} to $G_1\cup G'_2$ (with $\nu'$ playing the role of $\nu$) to obtain
a Hamilton decomposition of  $G_1\cup G'_2$. Thus ${\rm ham}(G)\ge f(n,\delta)-\eps n$.
This implies that ${\rm ham}^{\rm dir}(n,\delta)\ge f(n,\delta)-\eps n$.

The proof of~(ii) is identical to the proof of Theorem~\ref{mingraph}, except that instead of applying
Theorem~\ref{reggraph} we apply Theorem~\ref{regdigraph}.
\endproof

%%%%%%%%%%%%%%%%%%%%%%%%%%%%%%%%%%%%%%%%%%%%%%%%%%%%%%%%%%%%%%%%%%%%%%%%%%%%%%%%%%%%%

\section{Edge-disjoint Hamilton cycles in random tournaments, random and quasi-random graphs} \label{sec:tourn}

The purpose of this section is to prove Theorems~\ref{erdos},~\ref{Gnp} and~\ref{corquasi},
which give optimal packings of edge-disjoint Hamilton cycles in random tournaments, 
in very dense random graphs and in quasi-random graphs respectively.

\subsection{Optimal packings of Hamilton cycles in random tournaments}\label{sec:erdeos}
To prove Theorem~\ref{erdos}, we will apply the Max-Flow-Min-Cut theorem to show that a.a.s.~a random tournament~$T$
contains a $\delta^0(T)$-factor (see Lemma~\ref{tourn_reg}). Since a.a.s.~$\delta^0(T)$ is almost $|T|/2$, Theorem~\ref{erdos}
will then follow from Theorem~\ref{orientcor}. For our proof of Lemma~\ref{tourn_reg} we need some bounds on $\delta^0(T)$
as well as on the number of edges between any two sufficiently large subsets of vertices of $T$, which we will prove
in Lemma~\ref{tourn_edges}. To do this, we in turn need the following notation and some well-known facts about the
binomial distribution. Given $X \sim Bin(n-1, 1/2)$, we write
\begin{itemize}
\item $b(r) := \mathbb{P}\left(X = r \right) = {n-1 \choose r} (1/2)^{n-1}$,
\item $B(m_1, m_2) := \mathbb{P}\left(m_1 \leq X \leq m_2 \right)$, and
\item $B(m) := \mathbb{P}\left(X \leq m \right)$.
\end{itemize}
$b'(r)$, $B'(m_1, m_2)$ and $B'(m)$ are defined similarly for $X \sim Bin(n-2, 1/2)$.
\begin{prop} \label{binresults}
Suppose that $r \ge n/2 - \sqrt{2n \log n}$ and  $0 < h \le n^{3/5}$. Then
\begin{itemize}
\item[(i)] $\frac{b'(r)}{b(r)} \le 1+1/\log n$;
\item[(ii)] $b(n/2-h) \ge \frac{1}{2 \sqrt{n}} e^{-2h^2/n-4h^3/n^2}$;
\item[(iii)] $B(n/2-h) \le \frac{\sqrt{n}}{h} e^{-2h^2/n}$.
\end{itemize}
\end{prop}
(i) is a special case of Lemma~11(v) in~\cite{KnoxKO}. (ii) is a special case of Theorem~1.5 in~\cite{Bollobasbook}.
(iii) is a special case of the de Moivre-Laplace Theorem (Theorem~1.6(ii) in~\cite{Bollobasbook}).
The bounds in (ii) and (iii) are weaker than those in~\cite{Bollobasbook}, which allows us to replace appearances of the term $n-1$ (which would normally appear
e.g.~in the exponent of (iii)) by $n$.

\begin{lemma}\label{tourn_edges}
Suppose that $0 < \eps < 1$ is fixed.  
Let $T$ be a random tournament obtained from $K_n$ by orienting each
edge of $K_n$ randomly (where the probability for each of the two possible directions is $1/2$), independently of all other edges.
Then a.a.s.~the following conditions hold:
\begin{itemize}
\item[{\rm (i)}] $\delta^0(T)\ge n/2-(1+\eps)\sqrt{n(\log n)/2}$.
\item[{\rm (ii)}] $\delta^0(T)\le n/2-(1-\eps)\sqrt{n(\log n)/2}$.
\item[{\rm (iii)}] All vertex sets $A,B\subseteq V(T)$ satisfy $\left|e_T(A,B)-|A||B|/2\right|\le 9n^{3/2}$.
\end{itemize}
\end{lemma}
Assertions~(i) and~(ii) are very similar to those for the minimum degree of a binomial random graph $G_{n,1/2}$ on $n$ vertices with edge probability $1/2$. 
However, we cannot just quote these, as the correlations between vertex degrees in $G_{n,1/2}$ is slightly different than in a tournament.
For our purposes, it will suffice to apply Lemma~\ref{tourn_edges} with $\eps =1/10$, say.
\proof
First we prove (i).
Consider any vertex $x$. Let $h:=(1+\eps)\sqrt{n(\log n)/2}$. Then Proposition~\ref{binresults}(iii)
implies that
\begin{align*}
\pr(d^+_T(x) \le n/2-h)=B(n/2-h) & \le \frac{\sqrt{n}}{h} e^{-2h^2/n} \le n^{-(1+\eps)^2} \le n^{-1-2\eps}.
\end{align*}
The analogue holds for the indegree of $x$. So taking a union bound shows that the probability that $T$ does not satisfy (i) is at most
$2n\cdot n^{-1-2\eps} \le n^{-\eps}$.

Next we prove (ii). Let $h_1:=(1-\eps/2)\sqrt{n(\log n)/2}$ and $m_1:=n/2-h_1$.
Let $h_2:=(1-\eps)\sqrt{n(\log n)/2}$ and $m_2:=n/2-h_2$.
Let $Y$ be the number of vertices $x \in V(T)$ such that $m_1 \le  d^+_T(x) \leq m_2$.
Then Proposition~\ref{binresults}(ii) implies that 
\begin{align*}
\mathbb{E}(Y)  & \ge n (h_2- h_1)b(m_1) \ge
 n \cdot \frac{\eps}{2} \sqrt{n (\log n)/2} \cdot  \frac{1}{2 \sqrt{n}} e^{-2h_1^2/n-4h_1^3/n^2} \\
& \ge n \cdot (\log n)^{1/3} \cdot n^{-(1-\eps/2)^2} \ge n^{\eps/2}.
\end{align*}
Also
$$\mathbb{E}_2(Y) = n(n-1) \cdot 2 \cdot (1/2) \cdot B'(m_1-1, m_2-1)B'(m_1, m_2) \leq n^2 B'(m_1, m_2)^2.$$ 
Hence Proposition~\ref{binresults}(i) implies that 
\begin{equation}\label{eq:frac}
\frac{\sqrt{\mathbb{E}_2(Y)}}{\mathbb{E}(Y)} \leq \frac{\sum_{r = m_1}^{m_2} nb'(r)}{\sum_{r = m_1}^{m_2} nb(r)}
= \frac{\sum_{r = m_1}^{m_2} nb(r)\frac{b'(r)}{b(r)}}{\sum_{r = m_1}^{m_2} nb(r)} \leq 1 + \frac{1}{\log n}.
\end{equation}
So 
\begin{align}
Var(Y) &= \mathbb{E}_2(Y) + \mathbb{E}(Y) - \mathbb{E}(Y)^2 \leq (1 + 1/\log n)^2 \mathbb{E}(Y)^2 + \mathbb{E}(Y) - \mathbb{E}(Y)^2 \nonumber\\
&= (2/ \log n + 1/(\log n)^2)\mathbb{E}(Y)^2 + \mathbb{E}(Y)\label{eq:var}.
\end{align}
So by Chebyshev's inequality, 
$$\mathbb{P}(Y = 0) \leq \frac{Var(Y)}{\mathbb{E}(Y)^2} \leq \frac{2}{\log n} + \frac{1}{(\log n)^2} + \frac{1}{\mathbb{E}(Y)} \leq \frac{3}{\log n}.$$
Since the indegrees of the vertices satisfy the analogous properties, it follows that~(ii) fails with probability at most $6/\log n$.

Let us now check~(iii). We will first prove the following claim.

\smallskip

\noindent
\textbf{Claim.} \emph{Let $S,T\subseteq V(T)$ be such that $S\cap T=\emptyset$.
Then $$\pr(|e_T(S,T)-|S||T|/2|\ge 4 n^{3/2})\le 2e^{-2n}.$$}

\smallskip

\noindent To prove the claim, suppose first that least one of $S$, $T$ has size at most $4 n^{1/2}$. Then $e_T(S,T)\le 4n^{3/2}$ and so trivially
$\pr(|e_T(S,T)-|S||T|/2|\ge 4 n^{3/2})=0$. So we may assume that $|S|,|T|\ge 4n^{1/2}$. Let $a:=(12n/|S||T|)^{1/2}$.
Then $a<1$ and $a|S||T|/2\le 4n^{3/2}$. Thus Proposition~\ref{chernoff} implies that
\begin{align*}
\pr(|e_T(S,T)-|S||T|/2|\ge 4 n^{3/2})& \le \pr(|e_T(S,T)-|S||T|/2|\ge a|S||T|/2)\\
& \le 2e^{-a^2 |S||T|/6}= 2e^{-2n},
\end{align*}
which proves the claim.

\medskip

\noindent
Now consider any $A,B\subseteq V(T)$. Then
\begin{align*}
e_T(A,B) & =e_T(A,B\setminus A)+e_T(A\setminus B,A\cap B)+e_T(A\cap B,A\cap B)\\
& =e_T(A,B\setminus A)+e_T(A\setminus B,A\cap B)+\binom{|A\cap B|}{2}.
\end{align*}
Our claim implies that with probability at least $1-4e^{-2n}$
both $$|e_T(A,B\setminus A)-|A||B\setminus A|/2|\le  4 n^{3/2}$$ and
$$|e_T(A\setminus B,A\cap B)-|A\setminus B||A\cap B|/2|\le  4 n^{3/2}$$ hold.
But
\begin{align*}
\frac{|A||B\setminus A|}{2} & + \frac{|A\setminus B||A\cap B|}{2}+\binom{|A\cap B|}{2}+\frac{|A\cap B|}{2}=\\
%= \frac{1}{2}\left(|A||B\setminus A|+ |A\setminus B||A\cap B|+|A\cap B|^2\right)\\
& = \frac{1}{2}\left(|A||B\setminus A|+ (|A\setminus B|+|A\cap B|)|A\cap B|\right)\\
& = \frac{1}{2}\left(|A||B\setminus A|+ |A||A\cap B|\right)= \frac{|A||B|}{2}.
\end{align*}
Altogether this shows that with probability at least $1-4e^{-2n}$ we
have $$|e_T(A,B)-|A||B|/2|\le 8n^{3/2}+|A\cap B|/2\le 9n^{3/2}.$$
So taking a union bound shows that (iii) fails with probability at most $2^{2n}\cdot 4e^{-2n}$. 

Altogether this shows that with probability at most $n^{-\eps}+ 6/\log n+ 2^{2n}\cdot 4e^{-2n} \le 1/2$ at least
one of (i), (ii), (iii) fails.
\endproof 

\begin{lemma}\label{tourn_reg}
Suppose that $0<1/n\ll \eps\ll 1$. Let $T$ be a tournament on $n$ vertices which satisfies conditions (i)--(iii) of
Lemma~\ref{tourn_edges}. Then $T$ contains a $\delta^0(T)$-factor.
\end{lemma}
\proof
The proof is similar to that of Lemma~\ref{factor}(i).
Let $\delta:=\delta^0(T)$. As in the proof of Lemma~\ref{factor}(i), our aim is to apply the Max-Flow-Min-Cut theorem.
So let $H$ be the (unoriented) bipartite graph whose vertex classes $A$ and $B$ are both copies of $V(T)$ and in which
$a\in A$ is joined to $b\in B$ if $ab$ is a (directed) edge of $T$. Give every edge of $H$ capacity~1. Add a source $u^*$ which is
joined to every vertex $a\in A$ with an edge of capacity $\delta$. Add a sink $w^*$ which is
joined to every vertex $b\in B$ with an edge of capacity $\delta$. 
Note that an integer-valued $\delta n$-flow corresponds to the desired spanning subdigraph of $T$. Thus by the
Max-Flow-Min-Cut theorem it suffices to show that every cut has capacity at least $\delta n$.

So consider a minimal cut $\mathcal{C}$. Let $U$ be the set of all those vertices $a\in A$ for which $u^*a\notin \mathcal{C}$.
Similarly, let $W$ be the set of all those vertices $b\in B$ for which $bw^*\notin \mathcal{C}$.
Let $U':=A\setminus U$ and $W':=B\setminus W$.
Thus the capacity of $\mathcal{C}$ is
$$ c:=\delta |U'|+e_B(U,W)+ \delta |W'|.
$$
Note that $e_B(U,W)=e_T(U,W)\ge (\delta-|U'|) |W|$. Thus if $|W|\le \delta$, then
$$c=\delta |U'|+e_B(U,W)+ \delta |W'|\ge \delta |U'|-|U'||W|+\delta (|W|+|W'|)\ge \delta n,$$
as required.

The argument in the case when $|U|\le \delta$ is similar. So it remains to consider the case when $|U|,|W|\ge \delta$.
Note that Lemma~\ref{tourn_edges}(i),(ii) implies that 
\begin{equation} \label{deltamin}
n/2 - \sqrt{n \log n} \le \delta \le n/2 - \sqrt{n \log n}/2.
\end{equation}
In particular, this means
\begin{equation} \label{Ubound}
n/3 \le  \delta \le |U| \le n.
\end{equation}
So Lemma~\ref{tourn_edges}(iii) implies that
\begin{eqnarray}
c & \ge &\delta |U'|+\frac{|U||W|}{2}-9n^{3/2} + \delta |W'| \nonumber\\
& =& \delta |U'|+\frac{|U|}{2}(n-|W'|)-9n^{3/2} + \delta |W'|\nonumber \\
& \stackrel{(\ref{deltamin})}{\ge} & \delta |U'|+\frac{|U|}{2}(2\delta+\sqrt{n\log n}) -\frac{|U||W'|}{2} -9n^{3/2}+ \delta |W'|\nonumber\\
& = &\delta (|U'|+|U|)+\frac{|U|}{2}\sqrt{n\log n}-9n^{3/2} +|W'|\left(\delta-\frac{|U|}{2}\right)\nonumber\\
& \stackrel{(\ref{Ubound})}{\ge} & \delta n+\frac{n}{7}\sqrt{n\log n} +|W'|\left(\delta-\frac{|U|}{2}\right) \label{eq:cap1}\\
& \stackrel{(\ref{deltamin}),(\ref{Ubound})}{\ge} & \delta n+\frac{n}{7}\sqrt{n\log n} -|W'|\sqrt{n\log n}. \label{eq:cap2}
\end{eqnarray}
Now suppose first that $|W'|\le n-2\delta\le 2\sqrt{n\log n}$ (where the last inequality follows from (\ref{deltamin})).
Then (\ref{eq:cap2}) implies that
$$c\ge \delta n+\frac{n}{7}\sqrt{n\log n}-2n\log n\ge \delta n,$$
as required. So we are left with the case when $|U|,|W|\ge \delta$ and $|W'|\ge n-2\delta$. But then $|W|\le 2\delta$.
Moreover, interchanging the roles of $U$ and $W$ shows that instead of (\ref{eq:cap1}) one also has that
$$c\ge \delta n+\frac{n}{7}\sqrt{n\log n} +|U'|\left(\delta-\frac{|W|}{2}\right)\ge  \delta n+\frac{n}{7}\sqrt{n\log n}\ge \delta n,$$
as required.
\endproof

Theorem~\ref{erdos} is now an immediate consequence of Theorem~\ref{orientcor} and Lemmas~\ref{tourn_edges} and~\ref{tourn_reg}.

\removelastskip\penalty55\medskip\noindent{\bf Proof of Theorem~\ref{erdos}. }
Consider a random tournament $T$ obtained from $K_n$ by orienting each
edge of $K_n$ randomly (where the probability for each of the two possible directions is $1/2$), independently of all other edges.
Note that each tournament on $n$ vertices occurs equally likely in this model. So this is equivalent to choosing a tournament $T$ on
$n$ vertices uniformly at random.
Lemmas~\ref{tourn_edges} and~\ref{tourn_reg} together imply that a.a.s.~$T$ contains a $\delta^0(T)$-factor $T'$.
Since a.a.s.~$\delta^0(T)\ge n/2-\sqrt{n\log n}$ by Lemma~\ref{tourn_edges}(i), we can apply Theorem~\ref{orientcor} to obtain
a Hamilton decomposition of~$T'$. The Hamilton cycles in this Hamilton decomposition are as required in Theorem~\ref{erdos}.
\endproof

%%%%%%%%%%%%%%%%%%%%%%%%%%%%%%%%%%%%%%%%%%%%%%%%%%%%%%%%

\subsection{Optimal packings of Hamilton cycles in very dense random graphs}
Recall from Section~\ref{intro} that Theorem~\ref{Gnp} is only open in the case when $p$ tends to~$1$ rather quickly,
more precisely when $p\ge 1-n^{-1/4}(\log n)^9$. 
The proof of this case of Theorem~\ref{Gnp} is similar to that of Theorem~\ref{erdos}. The main step is to prove Lemma~\ref{Gnpcor}
below. Since a.a.s. $\delta(G_{n,p})\ge 7n/12$ (say) if $p\ge 2/3$, Theorem~\ref{reggraph} and
Lemma~\ref{Gnpcor} together immediately imply that Theorem~\ref{Gnp} holds if $p=p(n)\ge 2/3$.
It would not be difficult to replace the `2/3' by an arbitrary constant, though in this case one would have to apply 
Theorem~\ref{undir_decomp} rather than Theorem~\ref{reggraph}.
Below, we say  that a matching in a graph $G$ is \emph{optimal} if it covers all but at most one vertex of $G$.

\begin{lemma}\label{Gnpcor}
Suppose that $p=p(n)\ge 2/3$ is monotone. Then a.a.s. $G_{n,p}$ satisfies the following properties:
\begin{itemize}
\item[{\rm (i)}] If $\delta(G_{n,p})$ is even then $G_{n,p}$ contains a $\delta(G_{n,p})$-factor.
\item[{\rm (ii)}] If $\delta(G_{n,p})$ is odd then there is an optimal matching $M$ in $G_{n,p}$ such
that $G_{n,p}-M$ contains a $(\delta(G_{n,p})-1)$-factor.
\end{itemize}
\end{lemma}
Note that Lemma~\ref{Gnpcor} allows us to deduce not just Theorem~\ref{Gnp} (in the range when $p\ge 2/3$)
but the following stronger \emph{property} $\mathcal{H}$,
where a graph $G$ has property $\mathcal{H}$ if $G$ contains 
$\lfloor \delta(G)/2 \rfloor$ edge-disjoint Hamilton cycles, together with 
an additional edge-disjoint optimal matching if $\delta(G)$ is odd.
Property $\mathcal{H}$ was also verified for the relevant range of $p$ in~\cite{KnoxKO}, but not in~\cite{KrS}.

The next lemma guarantees a $\delta(G)$-factor or a $(\delta(G)-1)$-factor together with an optimal matching
in any graph $G$ whose minimum degree is close to $n$ and which satisfies some very weak conditions on the number of
vertices of degree $\delta(G)$ and on the size of the gap $\Delta(G)-\delta(G)$. 

\begin{lemma}\label{matchings}
Let $s,t,n\in\mathbb{N}$.
Suppose that $G$ is a graph on $n$ vertices such that $\delta(G)\ge n-t$ and let $s$ denote the number of vertices of degree $\delta(G)$
in $G$. Suppose that $n\ge s+3t+2t(\Delta(G)-\delta(G))$. Then the following properties hold.
\begin{itemize}
\item If $\delta(G)$ is even then $G$ contains a $\delta(G)$-factor.
\item If $\delta(G)$ is odd then there is an optimal matching $M$ in $G$ such that $G-M$ contains a $(\delta(G)-1)$-factor.
\end{itemize}
\end{lemma}
\proof
Let $\delta:=\delta(G)$. 
In our proof we will often use the following observation:

\textno
Whenever $G'$ is a spanning subgraph of $G$ with $\delta(G')=\delta$ and $X\subseteq V(G)$ is a set
of at least $2t$ vertices, then $G'[X]$ contains an optimal matching. Moreover, if $|X|$ is odd then any vertex
can be chosen as the one not covered by this matching. &(*)

(To see that $(*)$ holds, note that $\delta(G'[X])\ge |X|-t\ge |X|/2$. So $(*)$ follows from Dirac's theorem.)
Let $x_1,\dots,x_n$ be an enumeration of the vertices of $G$ such that $d(x_1)\ge \dots\ge d(x_n)$.
Let $X^{\rm min}$ be the set of all vertices of degree $\delta$ in $G$.
Let $d:=d(x_{2t})$. Note that $s < n-2t$, so $d>\delta$.
For each $i<2t$ we pick a set $N_i$ of $d(x_i)-d$ neighbours of $x_i$
in $V(G)\setminus (X^{\rm min}\cup \{x_1,\dots,x_{2t}\})$ such that these sets $N_i$ are pairwise disjoint.
To see that this can be done, suppose that we have already chosen $N_1,\dots,N_i$ for some $i<2t-1$ and that
we now wish to choose $N_{i+1}$. Since $X^{\rm min}\cup \{x_1,\dots,x_{2t}\}\cup N_1\cup\dots\cup N_i$
has size
$$s+2t+\sum_{j=1}^i|N_j|\le s+2t+(2t-2)(\Delta(G)-\delta)\le \delta-(\Delta(G)-\delta),$$
it follows that there are at least $\Delta(G)-\delta\ge d(x_i)-d$ possible vertices for $N_{i+1}$.
So we can choose $N_{i+1}$.
Let $G_0$ be the graph obtained from $G$ by deleting the edges between $x_i$ and $N_i$ for each $i<2t$.
So $\Delta(G_0)=d$ and $\delta(G_0)=\delta$. Let $X^{\rm max}_0$ be the set of all those vertices which have degree $d$ in $G_0$.
Thus all of $x_1,\dots,x_{2t}$ lie in $X^{\rm max}_0$ and so $|X^{\rm max}_0|\ge 2t$.
Let $X^{\rm min}_0$ be the set of all those vertices which have degree $\delta$ in $G_0$.
Thus $X^{\rm min}_0\subseteq X^{\rm min}\cup\bigcup_{i=1}^{2t-1}N_i$ and so
\begin{equation}\label{eq:s'}
|X^{\rm min}_0|\le s+\sum_{i=1}^{2t-1}|N_i|\le s+(2t-1)(\Delta(G)-\delta)=:s'.
\end{equation}
We will now prove the following claim.

\medskip

\noindent
\textbf{Claim.} \emph{Let $0\le i<d-\delta$. Suppose that $G_i$ is a spanning subgraph of $G$ with $\Delta(G_i)=d-i$ and $\delta(G_i)=\delta$.
Suppose that the set $X^{\rm max}_i$ of all vertices of degree $\Delta(G_i)$ in $G_i$ has size at least $2t$ and that the set
$X^{\rm min}_i$ of all vertices of degree $\delta$ in $G_i$ has size at most $s'+i$. Then the following holds:
\begin{itemize}
\item[{\rm (i)}] If $i\le d-\delta-2$ then $G_i$ contains a spanning subgraph $G_{i+1}$ with $\Delta(G_{i+1})=d-i-1$, $\delta(G_i)=\delta$
and such that the set $X^{\rm max}_{i+1}$ of all vertices of degree $\Delta(G_{i+1})$ in $G_{i+1}$ has size at least $2t$ and the set
$X^{\rm min}_{i+1}$ of all vertices of degree $\delta$ in $G_{i+1}$ has size at most $s'+i+1$.
\item[{\rm (ii)}] If $i=d-\delta-1$ and if $\delta$ is even then $G_i$ contains a $\delta$-factor.
\item[{\rm (iii)}] If $i=d-\delta-1$ and if $\delta$ is odd there is an optimal matching $M$ in $G_i$ such that $G_i-M$ contains a $(\delta-1)$-factor.
\end{itemize}}

\smallskip

Clearly, if we apply the claim $d-\delta$ times, starting  with $G_0$ then this gives a $\delta$-factor or an optimal matching and
a $(\delta-1)$-factor which are as required in the lemma.
So it suffices to prove the claim. First apply $(*)$ to find an optimal matching $M'$ in $G_i[X^{\rm max}_i]$. If
$|X^{\rm max}_i|$ is even, let $G_{i+1}$ be obtained from $G_i$ by deleting the edges in $M'$.
If $i\le d-\delta-2$ then $G_{i+1}$ satisfies~(i). If $i=d-\delta-1$ then $G_{i+1}$ is a $\delta$-factor.
So if $\delta$ is even then (ii) holds. Moreover, if $\delta$ is odd, this implies that $n$ must be even and we can apply $(*)$ to find a perfect
matching $M$ in $G_{i+1}$. Deleting $M$ gives a $(\delta-1)$-factor. So (iii) holds.

So we may assume that $|X^{\rm max}_i|$ is odd. Let $y\in X^{\rm max}_i$
be the vertex not covered by $M'$. If $i\le d-\delta-2$, choose a neighbour $y'$ of $y$ such that $y'\notin X^{\rm min}_i$ and
if $|X^{\rm max}_i|=2t$ then $y'\notin X^{\rm max}_i$.
This is possible since
$$\delta\ge 2t+s+2t(\Delta(G)- \delta)\stackrel{(\ref{eq:s'})}{=}2t+s'+\Delta(G)- \delta \ge 2t+(s'+d-\delta)>2t+|X^{\rm min}_i|.$$
Then the graph $G_{i+1}$ obtained from $G_i$ by deleting $M'$ and the edge $yy'$ satisfies~(i).
Note that if $i=d-\delta-1$ then $\delta$ must be odd (since then $|X^{\rm max}_i|$ is odd and $G_i$ has precisely $|X^{\rm max}_i|$
vertices of degree $\delta+1$ while all other vertices have degree $\delta$).
Choose any neighbour $y'$ of $y$. Let $G'_i$ be the graph obtained from $G_i$ by deleting $M'$ and $yy'$.
Then $y'$ has degree $\delta-1$ in $G'_i$ and all other vertices have degree $\delta$. But this implies that $n$ is odd
and so we can apply~$(*)$ to find a matching $M$ in $G'_i$ which covers all vertices apart from~$y'$.
Deleting $M$ from $G'_i$ yields a $(\delta-1)$-factor. So (iii) holds. 
\endproof

Similarly as at the beginning of Section~\ref{sec:erdeos}, given $X \sim Bin(n-1, p)$, we write
\begin{itemize}
\item $b(r) := \mathbb{P}\left(X = r \right) = {n-1 \choose r} p^r(1-p)^{n-1-r}$, and
\item $B(m) := \mathbb{P}\left(X \leq m \right)$.
\end{itemize}
$b'(r)$ and $B'(m)$ are defined analogously for $X \sim Bin(n-2, p)$.
The next lemma will be used to show that if $p\ge 1-n^{-2/3}$ then a.a.s.~$G_{n,p}$ satisfies the assumptions
of Lemma~\ref{matchings}.

\begin{lemma}\label{Gnpdegree}
Suppose that $0<p=p(n)\le n^{-2/3}$ is monotone. There exists an $\eps>0$ such that a.a.s $G_{n,p}$ satisfies the following properties:
\begin{itemize}
\item[{\rm (i)}] $\Delta(G_{n,p})\le 2n^{1/3}$.
\item[{\rm (ii)}] Either at most $(1-\eps)n$ vertices in $G_{n,p}$ have degree $\Delta(G_{n,p})$ or at least $\eps n$ vertices
of $G_{n,p}$ are isolated (or both).
\end{itemize}
\end{lemma}
Obviously (ii) is very crude, but we were not able to find an explicit statement in the literature which implies (ii), so we include a proof for completeness.
\proof
Condition~(i) is an extremely weak version of Corollary~3.4 in~\cite{Bollobasbook}. In order to check~(ii), let $q:=1-p$.
Suppose
first that $pqn\to \infty$ as $n\to \infty$. Let $Y$ denote the number of vertices of degree at most $pn$ in $G_{n,p}$.
Then the de Moivre-Laplace theorem implies that we have $\ex(Y)\in [n/3,2n/3]$. Note that
$$\ex_2(Y)=n(n-1)(p B'(np-1)^2+q B'(np)^2)\le n^2 B'(np)^2.$$
Moreover, a straightforward calculation shows that%
    \COMMENT{\begin{align*}\frac{b'(r)}{b(r)} &= \frac{{n-2 \choose r} p^{r} (1-p)^{n-r-2}}{{n-1 \choose r} p^r (1-p)^{n - r - 1}} = \frac{(n-2)!r! (n-r-1)!}{(n-1)!r!(n-r-2)!(1-p)}\\
& = \frac{n-r-1}{(n-1)(1-p)} = 1 + \frac{np - p - r}{(n-1)(1-p)}\le 1 + \frac{np}{(n-1)(1-p)}\le 1+\frac{3np}{2(n-1)}\le 1+2p.
\end{align*}
(Here the 2nd inequality holds since $p\le 1/3$.)}
$b'(r)/b(r)\le 1+2p$. Similarly as in (\ref{eq:frac}) one can use this to show that
$$\sqrt{\ex_2(Y)}/\ex(Y)\le 1+2p.$$
Similarly as in~(\ref{eq:var}) it follows that%
   \COMMENT{$Var(Y)= \mathbb{E}_2(Y) + \mathbb{E}(Y) - \mathbb{E}(Y)^2 \leq (1 + 2p)^2 \mathbb{E}(Y)^2 + \mathbb{E}(Y) - \mathbb{E}(Y)^2
\le 5p\ex(Y)^2+\ex(Y)$}
$Var(Y)\le 5p\ex(Y)^2+\ex(Y)$. %Choose $\eps$ and $\gamma$ such that $(1+\gamma)(1+\eps)/2 +\eps <1$.
Then Chebyshev's inequality implies that
$$\pr(Y\neq (1 \pm 1/6)\ex(Y))\le \frac{36 Var(Y)}{\ex(Y)^2}\le 36 \cdot 5p +\frac{36}{\ex(Y)}\to 0,$$
as $n\to \infty$. But if $Y= (1 \pm 1/6)\ex(Y)$, then $Y \in [5n/18,14n/18]$. This implies that there are at most
$n-Y\le 13n/18$ vertices of maximum degree.

So suppose next that $pqn\not\to \infty$. Let $Z$ denote the number of isolated vertices in $G_{n,p}$.
Thus there is a constant $\eps$ with $0 < \eps <1$ so that 
$$\ex(Z)=n(1-p)^{n-1}\ge ne^{-(p-p^2)n}\ge ne^{-pn/2}\ge 2\eps n.$$ 
Moreover,%
    \COMMENT{Note that $pqn\to \infty$ iff $pn\to \infty$ (since $q\ge 1/2$). But $p$ is monotone and so $pn$ is monotone.
So if $pqn\not\to \infty$ then there is some constant $C$ such that $pn\le C$.}
$$\ex_2(Z)=n(n-1)(1-p) b'(0)^2 \le n^2 b'(0)^2 \le n^2(1+2p)^2b(0)^2=(1+2p)^2\ex(Z)^2.$$
As before, we can use Chebyshev's inequality to show that the probability that at most $\eps n$ vertices are isolated
(i.e.~that $Z\le \eps n$) tends to~$0$, as $n\to \infty$.
\endproof

We can now combine Lemmas~\ref{matchings} and~\ref{Gnpdegree} to prove Lemma~\ref{Gnpcor}.

\removelastskip\penalty55\medskip\noindent{\bf Proof of Lemma~\ref{Gnpcor}. }
If $pn(1-p)/(\log n)^6\to \infty$, then Lemma~\ref{Gnpcor} is a special case of Lemma~21 in~\cite{KnoxKO}
(see the remarks after the proof of that lemma to see that its conditions are satisfied for such $p$).
Thus we may assume that $p\ge 1-n^{-2/3}$ (with room to spare). 
Let $\delta:=\delta(G_{n,p})$.
Applying Lemma~\ref{Gnpdegree} to the complement of $G_{n,p}$ shows that a.a.s. $\delta\ge n-1-2n^{1/3}\ge n-3n^{1/3}$
and there exists an $\eps>0$ such that at least one of the following conditions hold:
\begin{itemize}
\item At most $(1-\eps)n$ vertices of $G_{n,p}$ have degree $\delta$.
\item At least $\eps n$ vertices of $G_{n,p}$ have degree $n-1$.
\end{itemize}
Thus either $G_{n,p}$ is complete (in which case the lemma holds) or at most $(1-\eps)n=:s$ vertices have degree $\delta$.
Let $t:=3n^{1/3}$. Then
$$s+3t+2t(\Delta(G_{n,p})-\delta)\le (1-\eps)n+9n^{1/3}+6n^{1/3}(n-(n-3n^{1/3}))<n.$$
Thus Lemma~\ref{matchings} implies a.a.s.~$G_{n,p}$ satisfies~(i) or~(ii) of Lemma~\ref{Gnpcor}.
\endproof

%%%%%%%%%%%%%%%%%%%%%%%%%%%%%%%%%%%%%%%%%%%%%%%%%%%%%%%%%%%%%%%%%%%%%%%%%%%%%%

\subsection{Hamilton decompositions of quasi-random graphs}
Our aim of this section is to observe that if $\lambda$ is sufficiently small then every $(n,d,\lambda)$-graph
is a robust outexpander. We need the following notation.
Given sets $A,B$ of vertices of a graph $G$, we let $e'(A,B)$ denote the set of ordered pairs of vertices $a,b$ such that $a\in A$, $b\in B$
and $ab$ is an edge of~$G$. (So if $A=B$, we have $e'(A,A)=2e(A,A)$.) The following well-known result (see e.g.~Corollary~2.5 in Chapter 9 
of~\cite{ASp}) shows that for $(n,d,\lambda)$-graphs with small $\lambda$, $e'(A,B)$ is close to its `expected value'.

\begin{theorem} \label{quasibound}
Suppose that $G$ is an $(n,d,\lambda)$-graph.
Then for every pair $A,B$ of nonempty sets of vertices we have
$$
\left| \frac{e'(A,B)}{|A||B|} - \frac{d}{n} \right| \le \frac{\lambda}{\sqrt{|A||B|}}.
$$
\end{theorem}

\begin{lemma}\label{quasi_rob}
Suppose that $1/n \ll \eps \ll \nu \ll \tau ,\alpha<1$.
Suppose that $G$ is an $(n,d,\lambda)$-graph with $\lambda \le \eps n$ and $d\ge \alpha n$.
Then $G$ is a robust $(\nu,\tau)$-expander.
\end{lemma}
\proof 
Consider any set $S$ with $\tau n \le |S| \le (1- \tau) n$.
Let $RN =RN_\nu(S)$.
Then
$$e'(S,RN) \ge |S|d -\nu n^2 \ge (1-\tau/4)|S|d.$$ Since clearly $e'(S,RN)\le |S||RN|$ it follows
that $|RN|\ge d/2$. Thus $\lambda/\sqrt{|S||RN|}\le \sqrt{\eps}$.
Together with Theorem~\ref{quasibound} this implies that
$$
e'(S,RN) \le (d/n+\sqrt{\eps})|S||RN|  \le (d+ \sqrt{\eps} n)(1-\tau)|RN| \le (1-\tau/2)d|RN|.
$$
Altogether, this implies that $|RN| \ge (1+\tau/4)|S| \ge |S|+ \nu n$, as required.
\endproof

\removelastskip\penalty55\medskip\noindent{\bf Proof of Theorem~\ref{corquasi}. }
Let $\tau:=\tau(\alpha)$, where $\tau(\alpha)$ is as defined in Theorem~\ref{undir_decomp}.
Choose new constants $n_0\in\mathbb{N}$ and $\nu, \eps$ such that $0<1/n_0\ll \eps \ll \nu\ll \tau,\alpha$. 
Consider any $(n,d,\lambda)$-graph $G$ on $n\ge n_0$ vertices with $\lambda\le \eps n$ and such that $d\ge \alpha n$ is even.
Then Lemma~\ref{quasi_rob} implies that $G$ is a robust $(\nu,\tau)$-expander.
Thus Theorem~\ref{undir_decomp} implies that $G$ has a Hamilton decomposition.
\endproof

%%%%%%%%%%%%%%%%%%%%%%%%%%%%%%%%%%%%%%%%%%%%%%%%%%%%%%%%%%%%%%%%%%%%%%%%%%%%%%

\section{Algorithmic aspects and the Regularity lemma} \label{sec:algo}
 
In this section, we use Szemer\'edi's regularity lemma~\cite{reglem} to provide an algorithm which 
checks robust expansion and to give an algorithmic proof of Theorem~\ref{undir_decomp}.
 
\subsection{An algorithm for checking robust expansion} 
The purpose of this subsection is to prove Theorem~\ref{algo}, i.e.~our aim is to show that one can decide in polynomial time whether a digraph $G$ is a robust outexpander.
More precisely, there is a polynomial algorithm which either decides that $G$ is not a robust $(\nu,\tau)$-outexpander or
that $G$ is a robust outexpander with slightly worse parameters.

Roughly speaking, in order to prove Theorem~\ref{algo} we apply the regularity lemma to $G$ to obtain a reduced digraph $R$ and check whether $R$ is a robust outexpander.
The latter can be done in constant time. Lemma~\ref{robustG} shows that if $R$ is a robust outexpander, then
$G$ is a robust outexpander (with slightly worse parameters). Lemma~\ref{robustR} shows that if $G$ is a robust outexpander, then
$R$ is one as well.

Before we can state the regularity lemma for digraphs, we need the following notation.
If $G$ is an undirected bipartite graph with vertex classes $X$ and $Y$, then the
\emph{density} of $G$ is defined as
$$ d(X, Y) := \frac{e(X,Y)}{|X||Y|}.$$
Given $\eps >0$, we say that $G$ is \emph{$\eps$-regular} if for any $X'
\subseteq X$ and $Y' \subseteq Y$ with $|X'| \geq \eps |X|$ and $|Y'| \geq \eps
|Y|$ we have $|d(X',Y') - d(X, Y)| < \eps$. $G$ is \emph{$(\eps,d)$-regular} if $G$ is
$\eps$-regular and has density $d\pm \eps$.

Given disjoint vertex sets $X$ and $Y$ in a digraph $G$, we use $G[X,Y]$ to denote
the bipartite subdigraph of $G$ whose vertex classes are $X$ and $Y$ and whose
edges are all the edges of $G$ directed from $X$ to $Y$. We often view $G[X,Y]$ as
an undirected bipartite graph. In particular, we say $G[X,Y]$ is \emph{$\eps$-regular
with density $d$} if this holds when $G[X,Y]$ is viewed as an undirected graph.

Next we state
the degree form of the regularity lemma for (di-)graphs. A regularity lemma for
digraphs was proven by Alon and Shapira~\cite{AS}. The degree form follows from
this in the same way as the undirected version (see~\cite{BCCsurvey}
for a sketch of the latter). 
An algorithmic version of the (undirected) regularity lemma was proved in~\cite{Aetal}.
An algorithmic version of the directed version can be proved in essentially the same way
(see~\cite{CKKO} for a sketch of the argument proving a similar statement).

\begin{lemma}[Regularity lemma for (di-)graphs] \label{regularity_lemma}
For any $\eps, M'$ there exist $M, n_0$ such that if~$G$ is a digraph on $n \geq
n_0$ vertices and $d \in [0, 1]$, then there exists a partition of $V(G)$ into
$V_0, \dots, V_k$ and a spanning subdigraph $G'$ of $G$ such that the following conditions hold:
\begin{itemize}
\item[(i)] $M' \leq k \leq M$.
\item[(ii)] $|V_0| \leq \eps n$.
\item[(iii)] $|V_1| = \dots = |V_k| =: m$.
\item[(iv)] $d^+_{G'}(x) > d^+_{G}(x) - (d+\eps) n$ for all vertices $x \in V(G)$.
\item[(v)] $d^-_{G'}(x) > d^-_{G}(x) - (d+\eps) n$ for all vertices $x \in V(G)$.
\item[(vi)] For all $i \in [k]$ the digraph $G'[V_i]$ is empty.
\item[(vii)] For all $1 \leq i, j \leq k$ with $i \neq j$ the pair $G'[V_i,V_j]$ is
$\eps$-regular and either has density 0 or density at least~$d$. 
\end{itemize}
The analogue also holds when~$G$ is an undirected graph, in which case~(iv) and (v) are replaced by the condition that
$d_{G'}(x) > d_{G}(x) - (d+\eps) n$ for all vertices $x \in V(G)$.
\end{lemma}

Note that in the directed case the densities of the pairs $G'[V_i,V_j]$ and $G'[V_j,V_i]$ might be different from each other.
We refer to $V_0$ as the \emph{exceptional set} and to $V_1, \dots, V_k$ as \emph{clusters}. 
$V_0,V_1, \dots, V_k$ as above is also called a \emph{regularity partition} and $G'$ is
called the \emph{pure (di-)graph}.
Given a digraph $G$ on $n$ vertices, we form the \emph{reduced digraph $R$ of $G$
with parameters $\eps, d$ and $M'$} by applying the regularity lemma with these
parameters to obtain $V_0, \dots, V_k$. $R$ is then the digraph whose vertices are
the clusters $V_1,\dots,V_k$ and whose edges are those (ordered) pairs $V_iV_j$ of clusters
for which $G'[V_i,V_j]$ is nonempty. The \emph{reduced graph $R$} of an undirected graph $G$ is defined in a similar way.

The next result from~\cite{OS} implies that the property of a digraph~$G$ being a robust outexpander is `inherited'
by the reduced digraph $R$ of~$G$. It also shows that this even holds for spanning subdigraphs $R'$ of $R$ whose edges
correspond to all those pairs which have density a little larger than $d$. The latter statement will be used
in Section~\ref{sec:orient}. In~\cite{OS} the lemma is only stated in the directed case, but
the argument for the undirected case is identical. (A weaker version of this lemma was already proved
in~\cite{KOTchvatal}.)

\begin{lemma} \label{robustR}
Suppose that $1/ n \ll 1/M' \ll \eps \ll d \leq d^\prime \le \nu \le\tau < 1$ and $d'\le \nu/20$. Let $G$ be a digraph on $n$ vertices
which is a robust $(\nu , \tau)$-outexpander. Let $R$ be the reduced digraph of $G$ with parameters $\eps, d$ and $M'$
and let $G'$ be the pure digraph. Let $R'$ be the spanning subdigraph of $R$ such that $E(R')$ consists of all those edges $V_iV_j\in E(R)$ for which
$G'[V_i,V_j]$ has density at least~$d'$. Then $R'$ is a robust $(\nu/4 , 3 \tau)$-outexpander.
The analogue also holds if $G$ is an undirected robust $(\nu , \tau)$-expander.
\end{lemma}

The following lemma gives a converse to this (but with a weaker bound).

\begin{lemma}\label{robustG} 
Suppose that $1/n\ll \eps, 1/k \ll d,\nu,\tau < 1$. Let~$G$ be a
digraph on $n$ vertices. Let $V_0,V_1,\dots,V_k$
be a partition of $V(G)$ with $|V_0|\le \eps n$ and $|V_1|=\dots=|V_k|=:m$. Suppose that~$R$
is a digraph whose vertices are $V_1,\dots,V_k$ such that for every edge $V_iV_j$ of $R$ the
bipartite graph $G[V_i,V_j]$ is $\eps$-regular of density at least~$d$. If $R$ is a robust
$(\nu,\tau)$-outexpander, then $G$ is a robust $(d\nu^2/8,2\tau)$-outexpander.
Similarly, if $G$ is undirected and $R$ is a robust $(\nu,\tau)$-expander,
then $G$ is a robust $(d\nu^2/8,2\tau)$-expander.
\end{lemma}
\proof
We only prove the directed version, the argument for the undirected case is similar.
Consider any set $S \subseteq V(G)$ with $2\tau n \le |S| \le (1-2\tau)n$. 
Let $S_R$ be the set of all those clusters $V_i$ for which $|V_i \cap S| \ge \nu m/5$.
Note that
\begin{equation}\label{eq:SR}
|S_R| \ge (|S|-\nu m k/5 -|V_0|)/m \ge |S|/m -\nu k/4.
\end{equation}
So in particular $|S_R|\ge \tau k$. Let $RN(S_R):=RN^+_{\nu, R}(S_R)$.
Since $R$ is a robust $(\nu,\tau)$-outexpander, it follows that at least one of the following two properties hold:
\begin{itemize}
\item[(a)] $|RN(S_R)| \ge |S_R|+\nu k$.
\item[(b)] $|S_R|> (1-\tau)k$ and therefore $|RN(S_R)| \ge (1-\tau+\nu)k$.
\end{itemize}
(The second part of~(b) follows by considering the robust outneighbourhood of a subset of $S_R$ of size $(1-\tau)k$.)
Let $N_G$ denote the set of vertices of $G$ lying in clusters of $RN(S_R)$.
If (a) holds then
$$
|N_G| \ge |S_R|m+\nu km \stackrel{(\ref{eq:SR})}{\ge} |S|-\nu n/4 +3\nu n/4 = |S|+\nu n/2.
$$
If (b) holds then
$$
|N_G| \ge (1-\tau+\nu)km\ge (1-\tau)n\ge |S|+\nu n/2.
$$
Now let $N^*:=N_G \setminus RN^+_{d\nu^2/8, G}(S)$. We will show that $|N^*| \le \nu n/4$.
(Together with the bound on $|N_G|$ this then implies the lemma.) So suppose that $|N^*| \ge \nu n/4$.
Then there is a cluster $V_i$ which contains at least $\nu m/4$ vertices of $N^*$. (So $V_i\in RN(S_R)$.)
Let $V_i^*:=V_i \cap N^*$. Since no vertex in $V_i^*$ lies in the robust outneighbourhood $RN^+_{d\nu^2/8, G}(S)$ of $S$,
the number of edges from $S$ to $V_i^*$ in $G$ is at most $d\nu^2 n |V_i^*|/8$.
On the other hand, since $V_i\in RN(S_R)$, the number of edges in $R$ from $S_R$ to $V_i$
is at least $\nu k$. Let $S^*_R$ denote the initial vertices (i.e.~the initial clusters) of these edges.
Since $S^*_R\subseteq S_R$, every $V_j\in S^*_R$ satisfies $|V_j\cap S|\ge \nu m/5$. Together with the fact that the edge
$V_jV_i\in E(R)$ corresponds to an $(\eps,d)$-regular pair $G[V_i,V_j]$ in $G$ and $|V_i^*| \ge \nu m/4\ge \eps m$, it follows that 
the number of edges in $G$ from $V_j\cap S$ to $V^*_i$ is at least
$5d|V_j\cap S||V_i^*|/6\ge d\nu m|V_i^*|/6$. Summing over all clusters $V_j\in S_R^*$ shows that the number of
edges in $G$ from $S$ to $V_i^*$ is at least $\nu k\cdot  d\nu m|V_i^*|/6\ge d\nu^2n|V_i^*|/7$, a contradiction.
\endproof

\removelastskip\penalty55\medskip\noindent{\bf Proof of Theorem~\ref{algo}. } Let $d:=\nu/20$.
Choose $\eps>0$ and $n_0,M'\in\mathbb{N}$ such that $1/n_0\ll \eps,1/M'\ll d$.
We may assume that $n\ge n_0$ (since otherwise we can solve the problem by complete enumeration).
Apply the regularity lemma (Lemma~\ref{regularity_lemma}) with parameters $\eps,d,M'$ to $G$ obtain a reduced digraph $R$
and a pure digraph $G'$.
Now check whether $R$ is a robust $(\nu/4,3\tau)$-outexpander (this can be done in constant time). If it is not
a robust $(\nu/4,3\tau)$-outexpander, then Lemma~\ref{robustR} (with $R$ and $d$ playing the roles of $R'$ and $d'$) implies that $G$ is not a robust $(\nu,\tau)$-outexpander.
If $R$ is a robust $(\nu/4,3\tau)$-outexpander, then Lemma~\ref{robustG} (with $G'$ playing the role of $G$) implies that $G'$ is
a robust $(\nu^3/2560,6\tau)$-outexpander. Since $G'\subseteq G$, this shows that $G$ is a robust $(\nu^3/2560,6\tau)$-outexpander.
The algorithm for (undirected) robust expanders is similar.
\endproof 

\subsection{Regular orientations -- algorithmic proof of Lemma~\ref{reggraphorient}} \label{sec:orient}

The purpose of this subsection is to give an algorithmic proof of the following version of Lemma~\ref{reggraphorient}.

\begin{lemma}\label{reggraphorient2}
Suppose that $0<1/n\ll \nu'\ll \nu\le \tau\ll \alpha<1$.
Let $G$ be an $r$-regular graph on $n$ vertices such that $r\ge \alpha n$ is even and $G$ is a robust $(\nu,\tau)$-expander.
Then one can orient the edges of $G$ in such a way that the oriented graph $G^{\rm orient}$ thus obtained from $G$ is an
$r/2$-regular robust $(\nu',6\tau)$-outexpander. Moreover, this orientation $G^{\rm orient}$ can be found in time polynomial in~$n$.
\end{lemma}
Note that Lemma~\ref{reggraphorient} guarantees that $G^{\rm orient}$ is robust $(\nu',\tau)$-outexpander, instead of
just a robust $(\nu',6\tau)$-outexpander. But this does not matter for the applications. In particular, using a similar argument as
before, Lemma~\ref{reggraphorient2} immediately gives an algorithmic proof of Theorem~\ref{undir_decomp}, as the only other ingredient is the 
(algorithmic) Theorem~\ref{decomp}.

The proof idea in this case is to consider a suitable `fractional' and `almost regular' orientation of the reduced graph~$R$ of~$G$ obtained from the
regularity lemma.
This orientation is `lifted' to~$G$ and minor `irregularities' are eliminated by re-orienting a small proportion of the edges
(via Lemma~\ref{pathswitch}). 

We will need the following observations, which follow immediately from the definition of robust outexpansion.

\begin{prop} \label{basics}
Suppose that $0 < 1/n \ll \beta \ll \nu \le \tau \ll \alpha <1$ and that 
$G$ is a robust $(\nu, \tau)$-outexpander on $n$ vertices.
\begin{itemize}
\item[(i)] 
 Suppose that $G'$ is obtained from $G$ by reversing the orientations of at most $\beta n^2$ edges.
 Then $G'$ is a robust $(\nu/2, \tau)$-outexpander.
\item[(ii)] Suppose that $\delta^0(G) \ge \alpha n$. Then $G$ has diameter at most $1/\nu$.
\end{itemize}
\end{prop}

For a vertex $x$ in an oriented graph $G$, let ${\rm disc}(x,G):=|d^+_G(x)-d^-_G(x)|$.
For an oriented graph $G$, let ${\rm disc}(G):=\sum_{x \in V(G)}{\rm disc}(x,G)$.
So~$G$ is regular if and only if ${\rm disc}(G)=0$.

\begin{lemma} \label{pathswitch}
Suppose that $0< 1/n \ll \beta \ll \nu\le \tau \ll \alpha < 1$ and that $\alpha n\in\mathbb{N}$. Let $G$ be
a $2\alpha n$-regular graph on $n$ vertices which has an orientation $G'$ which is a robust $(\nu, \tau)$-outexpander
with ${\rm disc}(G') \le \beta n^2$. Then $G$ has an orientation $G''$ 
which is an $\alpha n$-regular robust $(\nu/2, \tau)$-outexpander.
Moreover, this orientation can be found in time polynomial in~$n$.
\end{lemma}
The idea of the proof is to repeatedly apply the following argument:
as long as there is a vertex $x$ whose outdegree is less than $\alpha n$, then there must also be a vertex
$y$ whose indegree is greater than  $\alpha n$. We find a short directed path~$P$ from $x$ to~$y$ and reverse the orientation
of the edges on $P$.

\proof 
We say that a vertex $x$ in an oriented graph $H$ is \emph{$H$-balanced} if ${\rm disc}(x,H) \le \alpha n$.
Let $B$ be the set of $G'$-balanced vertices.
Note that ${\rm disc}(G') \le \beta n^2$ implies that 
\begin{equation} \label{Bee}
|B| \ge (1-\sqrt{\beta}) n.
\end{equation} 
Let $G_0:=G'$ and consider the following statements, where $0 \le i \le \beta n^2/2$:
\begin{itemize}
\item[(a$_i$)] If $i \ge 1$, then $G_i$ is an oriented graph obtained from $G_{i-1}$ by reversing the orientations of the edges
of a directed path $P_{i-1}$ of length at most $4/\nu$ in~$G_{i-1}$.
\item[(b$_i$)] ${\rm disc}(x,G_{i}) \le {\rm disc}(x,G_{0})$ for all $x\in V(G)$. 
\item[(c$_i$)] ${\rm disc}(G_i) \le \beta n^2-2i$.
\end{itemize}
Note that (a$_0$),(b$_0$) vacuously hold and (c$_0$) holds by assumption.

Suppose inductively that for some $i\ge 0$ we have constructed $G_0,\dots, G_i$ such that (a$_i$)--(c$_i$) hold.
Suppose also that there is a vertex $x$ with ${\rm disc}(x,G_i)>0$. Without loss of generality, we have
$d^+_{G_i}(x)> d^-_{G_i}(x)$. So $d^+_{G_i}(x)> \alpha n$ and thus
$x$ has an outneighbour $x' \in B$ (in~$G_i$). If $x \in B$, let $x^*:=x$, otherwise, let $x^*:=x'$.
There must also be a vertex $y$ with $d^-_{G_i}(y)> d^+_{G_i}(y)$.
So $d^-_{G_i}(y)> \alpha n$ and $y$ has an inneighbour $y' \in B$ (in~$G_i$).
If $y \in B$, let $y^*:=y$, otherwise, let $y^*:=y'$.

Note that (b$_i$) implies that any $b \in B$ is $G_i$-balanced. 
So it follows from~(\ref{Bee}) that $\delta^0(G_i[B]) \ge \alpha n/3 \ge \alpha |B|/3$.
Moreover,~(\ref{Bee}) implies that $G_i[B]$ is a robust $(\nu/2, 2\tau)$-outexpander.
So by Proposition~\ref{basics}(ii), $G_i[B]$ 
contains a (directed) $x^*y^*$-path $P'_i$ of length at most $2/\nu$. 
We now extend $P'_i$ into an $xy$-path~$P_i$ in~$G_i$ of length at most $2/\nu+2 \le 4/\nu$
(where we may have $P'_i=P_i$).
Let $G_{i+1}$ be obtained from $G_i$ by reversing the orientations of the edges along $P_i$.
Then ${\rm disc}(x,G_{i+1}) < {\rm disc}(x,G_i)$, ${\rm disc}(y,G_{i+1}) < {\rm disc}(y,G_i)$
and ${\rm disc}(z,G_{i+1})={\rm disc}(z,G_i)\le  {\rm disc}(z,G_0) $ for any vertex $z\neq x,y$.
So (a$_{i+1}$)--(c$_{i+1}$) hold.

So we must have ${\rm disc}(G_s)=0$ for some $s \le  \beta n^2/2$.
Let $G'':=G_s$. By definition, $G''$ is $\alpha n$-regular.
Moreover, $G''$ is obtained from $G'$ by reversing the orientations of at most $2\beta n^2/\nu\le \sqrt{\beta}n^2$ edges.
Proposition~\ref{basics}(i) now implies that $G'$ is a  robust $(\nu/2, \tau)$-outexpander, as required.
Moreover, each $P_i$ (and thus $G''$) can be found in polynomial time, since the shortest path problem can be solved in polynomial time.
\endproof 

The following lemma follows immediately from Lemma~4.10(iii) in~\cite{monster}. As remarked after its proof in~\cite{monster},
the proof of this lemma is algorithmic.%
   \COMMENT{In~\cite{monster} $(\eps^{1/12},d')$-regular means multiplicatively regular, but this also implies $(\eps^{1/12},d')$-regular
in our (additive) sense.} 
\begin{lemma}\label{randomregslice}
Suppose that $0<1/m\ll \eps\ll d'\le d\le 1$. Let $G$ be an $\eps$-regular bipartite graph of density $d$ with vertex classes of size $m$.
Then $G$ contains an $(\eps^{1/12},d')$-regular spanning subgraph $J$.
Moreover, if $x \in V(G)$ satisfies $d_G(x) =(d \pm \eps)m$, then $d_{J}(x) =(d' \pm \sqrt{\eps})m$.
\end{lemma}

\removelastskip\penalty55\medskip\noindent{\bf Proof of Lemma~\ref{reggraphorient2}. }
Choose $\eps,d,d_1>0$ and $M'\in\mathbb{N}$ such that $1/n\ll 1/M'\ll \eps \ll d \ll \nu'\ll d_1\ll \nu$.
Apply the undirected version of the regularity lemma (Lemma~\ref{regularity_lemma}) with parameters $\eps,d,M'$ to $G$ obtain a reduced graph $R$
and a pure graph $G'$. For  each edge $XY$ of $R$, let $d_{XY}$ denote the density of $G'[X,Y]$ and let $d'_{XY}:=d_{XY}/2$.
Let $R^*$ be the digraph obtained from $R$ by replacing each (undirected) edge $XY$ of $R$ by 
the directed edges $XY$ and $YX$. Apply Lemma~\ref{randomregslice} with $G'[X,Y]$, $d_{XY}$ and  $d'_{XY}$ playing the roles $G$, $d$ and $d'$ 
to obtain an $(\eps^{1/12},d'_{XY})$-regular spanning subgraph $J_{XY}$ of $G'[X,Y]$.
We say that $x \in X \cup Y$ is \emph{$XY$-useful} if $d_{G'[XY]}(x) =(d_{XY} \pm \eps)m$.
So if $x$ is $XY$-useful, then by Lemma~\ref{randomregslice} we have $d_{J_{XY}}(x) =(d'_{XY} \pm \sqrt{\eps})m$.

Let $J_{YX}:=G'[X,Y]\setminus E(J_{XY})$. Note that the above implies that $J_{YX}$ is  $(2\eps^{1/12},d'_{XY})$-regular
and that if $x$ is $XY$-useful, then $d_{J_{YX}}(x) =(d'_{XY} \pm 2\sqrt{\eps})m$.
For each edge $XY$ of $R$, orient all edges in $J_{XY}$ from $X$ to $Y$ and orient all edges in $J_{YX}$ from 
$Y$ to $X$. 

We say that an edge of $G$ is \emph{good} if it lies in $G'[X,Y]$ for some edge $XY$~of $R$
and \emph{bad} otherwise.
The above gives an orientation of all good edges of $G$.
Orient the bad edges of $G$ arbitrarily to obtain an oriented graph $G^*$.
Let $R_1$ be the spanning subgraph of $R$ such that $E(R_1)$ consists of all those edges $XY\in E(R)$ for which
$d_{XY}\ge d_1$. By the undirected version of Lemma~\ref{robustR}, $R_1$ is a robust $(\nu/4,3\tau)$-expander.
Let $R^*_1$ be the digraph obtained from $R_1$ by replacing each (undirected) edge $XY$ of $R_1$ by 
the directed edges $XY$ and $YX$. Clearly, it follows that $R^*_1$ is a robust $(\nu/4,3\tau)$-outexpander.
Moreover, note that for any edge $XY$ of $R^*_1$, the density of the corresponding pair $G^*[X,Y]$ is at least $d_1/3$.
Thus we can apply the directed version of Lemma~\ref{robustG} with $G^*,R^*_1,d_1/3$ playing the roles of $G,R,d$ to 
see that $G^*$ is a  robust $(d_1\nu^2/(3 \cdot 128),6\tau)$-outexpander and thus also a robust $(2\nu',6\tau)$-outexpander. 

We claim that  ${\rm disc}(G^*) \le 3dn^2$.
To prove the claim, for any vertex $x \in V(G)$, let $b(x)$ denote the number of bad edges
incident to $x$. Note that Lemma~\ref{regularity_lemma} implies that
\begin{equation} \label{bx}
\sum_{x \in V(G)} b(x) \le n(\eps+d)n+|V_0|n \le 2d n^2.
\end{equation}
Let $G^*_{\rm good}$ be the oriented subgraph of $G^*$ consisting of (orientations of) the good edges.
Suppose that $x \in X$ for some cluster $X$ and that $x$ is $XY$-useful. Then our construction implies that
\begin{equation} \label{good}
\left| |N^+_{G^*_{\rm good}}(x) \cap Y| -|N^-_{G^*_{\rm good}}(x) \cap Y \right| \le 3\sqrt{\eps} m.
\end{equation}
Note that for a given edge $XY$ of~$R$, the number of vertices $x \in X$ which are not $XY$-useful is at most $2\eps m$.
We say that a vertex $x \in X$ is \emph{useful} if it is $XY$-useful for at least at least $(1-\sqrt{\eps})|R|$ clusters $Y$.
Suppose that $x$ is useful. Then~(\ref{good}) implies that 
\begin{equation} \label{disc}
{\rm disc}(x,G^*) \le b(x)+3\sqrt{\eps} m|R| +\sqrt{\eps} m|R| \le b(x)+4\sqrt{\eps}n.
\end{equation}
Note that the total number of vertices which are not useful is at most%
    \COMMENT{Fix a cluster $X$. Count the number $p$ of all pairs $(x,Y)$ such that $x\in X$ and $x$ is not $XY$ useful.
Then $p\le 2\eps m|R|$. So $X$ contains at most $p/\sqrt{\eps}|R|\le 2\sqrt{\eps}m$ useless vertices.}
$|V_0|+2\sqrt{\eps}|R|m \le 3\sqrt{\eps} n$.
Also, if $x$ is not useful, we have ${\rm disc}(x,G^*) \le n$.
So~(\ref{bx}) and~(\ref{disc}) together imply that
${\rm disc}(G^*) \le 2d n^2+ 4\sqrt{\eps} n^2+ 3\sqrt{\eps}n^2 \le 3dn^2$.
So we can now apply Lemma~\ref{pathswitch} with $G^*,3d,2\nu',6\tau, r/2n$ playing the roles of $G',\beta,\nu,\tau, \alpha$
to obtain an $r/2$-regular orientation of $G^{\rm orient}$ of $G$ which is a robust $(\nu',6\tau)$-outexpander.
\endproof

\section{Acknowledgements}

We are grateful to the referee for suggesting a way of `derandomizing' the proof of Lemma~\ref{reggraphorient}
(and thus of Theorem~\ref{undir_decomp}).

\medskip

{\footnotesize \obeylines \parindent=0pt

Daniela K\"{u}hn, Deryk Osthus 
School of Mathematics
University of Birmingham
Edgbaston
Birmingham
B15 2TT
UK
}
\begin{flushleft}
{\it{E-mail addresses}:
\tt{\{d.kuhn,d.osthus\}@bham.ac.uk}}
\end{flushleft}

\end{document}